\documentclass[preprint,3p,12pt]{elsarticle}
\usepackage{amsmath}
\usepackage{stmaryrd}
\usepackage{bbding}
\usepackage{dcolumn}
\usepackage{graphicx}
\usepackage{amsfonts}
\usepackage{amssymb}
\usepackage{psfrag}
\usepackage{wrapfig}
\usepackage{subfigure}
\usepackage{makeidx}
\usepackage{bm}
\usepackage{epsf}
\usepackage{epsfig}
\usepackage{setspace}
\usepackage{color}

\newtheorem{rmk}{Remark}

\begin{document}
\title{Three-dimensional discontinuous Galerkin based high-order gas-kinetic scheme and GPU implementation}
\author[BNU]{Yuhang Wang}
\ead{hskwyh@outlook.com}
\author[BNU]{Liang Pan\corref{cor}}
\ead{panliang@bnu.edu.cn}
\address[BNU]{Laboratory of Mathematics and Complex Systems, School of Mathematical Sciences, Beijing Normal University, Beijing, China}
\cortext[cor]{Corresponding author}

\begin{abstract}
In this paper, the discontinuous Galerkin based high-order
gas-kinetic schemes (DG-HGKS) are developed for the
three-dimensional Euler and Navier-Stokes equations. Different from
the traditional discontinuous Galerkin (DG) methods with Riemann
solvers, the current method adopts a kinetic evolution process,
which is provided by the integral solution of Bhatnagar-Gross-Krook
(BGK) model. In the weak formulation of DG method, a time-dependent
evolution function is provided, and both inviscid and viscous fluxes
can be calculated uniformly. The temporal accuracy is achieved by
the two-stage fourth-order discretization, and the second-order
gas-kinetic solver is adopted for the fluxes over the cell interface
and the fluxes inside a cell. Numerical examples, including accuracy
tests and Taylor-Green vortex problem, are presented to validate the
efficiency and accuracy of DG-HGKS. Both optimal convergence and
super-convergence are achieved by the current scheme. The comparison
between DG-HGKS and high-order gas-kinetic scheme with weighted
essential non-oscillatory reconstruction (WENO-HGKS) is also given,
and the numerical performances are comparable with the approximate
number of degree of freedom. To accelerate the computation, the
DG-HGKS is implemented with the graphics processing unit (GPU) using
compute unified device architecture (CUDA). The obtained results are
also compared with those calculated by the central processing units (CPU)
code in terms of the computational efficiency. The speedup of GPU
code suggests the potential of high-order gas-kinetic schemes for the
large scale computation.
\end{abstract}

\begin{keyword}
High-order gas-kinetic scheme (HGKS), discontinuous Galerkin (DG) method,
graphics processing unit (GPU).
\end{keyword}

\maketitle
\section{Introduction}
The numerical simulation of compressible flows is an important issue
for the computational fluid dynamics, and the construction of
high-order numerical schemes becomes demanding for the flows with
complicated structures, such as shock, vortex, boundary layer and
turbulence. In recent decades, there have been continuous interests
and efforts on the development of high-order schemes, including
discontinuous Galerkin (DG), spectral volume (SV) \cite{SV},
spectral difference (SD) \cite{SD}, flux reconstruction (FR)
\cite{FR}, correction procedure using reconstruction (CPR)
\cite{CPR-1,CPR-2}, essential non-oscillatory (ENO) \cite{ENO-1,ENO-2} and
weighted essential non-oscillatory (WENO) schemes
\cite{WENO-Liu,WENO-JS,WENO-Z}, etc. The DG method was originally
introduced by Reed and Hill \cite{DG0} in the framework of neutron
transport. Cockburn and Shu further developed the DG method in a
series of papers \cite{DG1,DG2,DG3}, in which a framework  was
established to solve the nonlinear time-dependent hyperbolic
conservation laws. Inspired by the great success of DG method for
the first-order system, a natural choice to solve a second-order
system is to convert it into a first-order system by introducing
additional variables, and to apply DG method to the first-order
system directly. There are mainly two kinds of approaches, i.e., the
first Bassi-Rebay (BR1) scheme \cite{BR1} and local discontinuous
Galerkin (LDG) method \cite{LDG}. To achieve the compactness,  the
second Bassi-Rebay (BR2) scheme \cite{BR2}  was developed based on
BR1 scheme. Recently, based on the direct weak formulation for
solutions of parabolic equations, a direct discontinuous Galerkin
(DDG) method was introduced to solve the diffusion equation and
Navier-Stokes equations \cite{DDG1,DDG2}. Generally, for the finite
volume methods, such as ENO and WENO schemes, a large stencil is
often used to achieve high-order accuracy. Due to the selection of
large candidate stencils, these methods face problems on the
unstructured meshes, especially for the three-dimensional
computation. Meanwhile, DG methods achieve the high-order accuracy
by reconstructing a piecewise discontinuous polynomial within each
cell, and only the immediate neighbors are involved for the update
of flow variables. Compared with the WENO-type methods, the DG
methods are less robust. The codes eventually blow up due to the
spurious oscillation near discontinuities \cite{DG2,DG3}. Many
limiters and trouble cell indicators have been introduced to
suppress the oscillations in the literatures for DG methods
\cite{DG-limiter-1,DG-limiter-2}.

In the past decades, the gas-kinetic schemes (GKS) have been
developed systematically based on the Bhatnagar-Gross-Krook (BGK)
model \cite{BGK-1,BGK-2} in the finite volume framework, and applied
successfully from low speed flows to hypersonic ones
\cite{GKS-Xu1,GKS-Xu2}. The gas-kinetic scheme presents a gas
evolution process from kinetic scale to hydrodynamic scale, where
both inviscid and viscous fluxes are recovered from a time-dependent
and genuinely multi-dimensional gas distribution function at cell
interface. Based on the two-stage fourth-order temporal
discretization for Lax-Wendroff type solvers
\cite{GRP-high-1,GRP-high-2}, the high-order gas-kinetic schemes
(HGKS) have been constructed and applied for the compressible flow
simulations \cite{GKS-high-1,GKS-high-2,GKS-high-3}. The
fourth-order and even higher-order temporal accuracy can be achieved
with the implementation of the traditional second-order or
third-order GKS evolution model. More importantly, HGKS is as robust
as the second-order scheme, and works perfectly from the subsonic to
hypersonic viscous heat conducting flows. With the implementation of
three-dimensional WENO reconstruction, the two-stage fourth-order
gas-kinetic scheme has been successfully implemented in the direct
numerical simulation (DNS) for compressible turbulences
\cite{GKS-high-4}. To improve the efficiency, a parallel code of
HGKS is developed \cite{GKS-high-4}, where the two-dimensional
domain decomposition and message passing interface (MPI) are used
for the implementation of parallel computing.  The scalability of
MPI code is examined up to $1024$ cores on the TianHe-II
supercomputer, and the MPI code scales properly with the number of
processors due to the explicit formulation of the algorithm. With
the parallel code, the HGKS provides us a powerful tool for the DNS
study from the subsonic to supersonic turbulences.

In this paper, a discontinuous Galerkin based high-order gas-kinetic
scheme (DG-HGKS) is developed for the compressible Euler and
Navier-Stokes equations. Different from the traditional DG methods
with Riemann solvers \cite{Riemann-appro}, a time-dependent
numerical flux is provided by the gas-kinetic flow solver. The
DG-HGKS with third-order and fourth-order spatial accuracy is
developed with $\mathbb P_2$ and $\mathbb P_3$ elements, where
$\mathbb P_k$ is the space of polynomials with the sum of degree at
most $k$ for the sum of all variables. Meanwhile, the temporal
accuracy is achieved by the two-stage fourth-order discretization
\cite{GRP-high-1,GKS-high-1}. In the computations, we mainly
concentrate on the smooth flows, and numerical examples are
presented to validate the performance of current method. The
accuracy tests are presented to validate the optimal convergence and
super-convergence for DG scheme with uniform and nonuniform meshes.
The Taylor-Green vortex problem is provided to validate the
efficiency and accuracy for the nearly incompressible turbulence. As
comparison, the numerical results of HGKS with weighted essential
non-oscillatory reconstruction (WENO-HGKS) are given as well. As
expected, the numerical performances, including the resolution and
efficiency, are comparable with the approximate number of degree of
freedom for two methods. To accelerate the computation, the DG-HGKS
is implemented to run on graphics processing unit (GPU) using
compute unified device architecture (CUDA). The computational
efficiency using single Nvidia TITAN RTX and Tesla V100 GPU is
demonstrated. Obtained results are compared with those obtained by
an octa-core Intel i7-9700K CPU in terms of calculation time.
Compared with the CPU code, 6x-7x speedup is achieved for TITAN RTX
and 10x-11x speedup is achieved for Tesla V100. Due to the explicit
formulation of HGKS, the GPU code can be implemented for WENO-HGKS
\cite{GKS-high-4} as well.  In the future, more challenging
compressible flow problems, such as the turbulent channel flows and
the flat plate turbulent boundary layer, will be investigated with
multiple GPUs.

This paper is organized as follows. In Section 2, the BGK equation
and GKS is briefly reviewed. The DG-HGKS is presented in
Section 3, and the GPU architecture and code design are given in
Section 4. Numerical examples are included in Section 5. The last
section is the conclusion.

\section{BGK equation and gas-kinetic scheme}
The three-dimensional BGK equation \cite{BGK-1,BGK-2} can be
written as
\begin{equation}\label{bgk-3d}
f_t+uf_x+vf_y+wf_z=\frac{g-f}{\tau},
\end{equation}
where $\boldsymbol{u}=(u,v,w)$ is the particle velocity, $f$ is the
gas distribution function, $g$ is the three-dimensional Maxwellian
distribution and $\tau$ is the collision time. The collision term
satisfies the compatibility condition
\begin{equation}\label{compatibility}
\int \frac{g-f}{\tau}\psi \text{d}\Xi=0,
\end{equation}
where
$\displaystyle\psi=(1,u,v,w,\frac{1}{2}(u^2+v^2+w^2+\xi^2))^T$, the
internal variables $\xi^2=\xi^2_1+...+\xi^2_K$,
$\text{d}\Xi=\text{d}u\text{d}v\text{d}w\text{d}\xi_1…\text{d}\xi_{K}$,
$\gamma$ is the specific heat ratio and  $K=(5-3\gamma)/(\gamma-1)$
is the degrees of freedom for three-dimensional flows. Taking
moments of BGK equation Eq.\eqref{bgk-3d}, the macroscopic equation
can be written as
\begin{align}\label{hyper}
\frac{\partial Q}{\partial t}+\nabla\cdot \boldsymbol{F}(Q)=0,
\end{align}
where the conservative variables $Q=(\rho, \rho U,\rho V, \rho W,
\rho E)^T$ and the fluxes $\boldsymbol{F}(Q)$ are given by
\begin{align*}
Q&=\displaystyle \int \psi f(\boldsymbol{x},t,\boldsymbol{u},\xi) \text{d}\Xi,\\
\boldsymbol{F}(Q)&=\displaystyle \int  \boldsymbol{u}\cdot
\boldsymbol{n} \psi f(\boldsymbol{x},t,\boldsymbol{u},\xi)
\text{d}\Xi.
\end{align*}
According to the Chapman-Enskog expansion for BGK equation, the
macroscopic governing equations can be derived. In the continuum
region, the BGK equation can be rearranged and the gas distribution
function can be expanded as
\begin{align*}
f=g-\tau D_{\boldsymbol{u}}g+\tau D_{\boldsymbol{u}}(\tau
    D_{\boldsymbol{u}})g-\tau D_{\boldsymbol{u}}[\tau
            D_{\boldsymbol{u}}(\tau D_{\boldsymbol{u}})g]+...,
\end{align*}
where $D_{\boldsymbol{u}}=\displaystyle\frac{\partial}{\partial
t}+\boldsymbol{u}\cdot \nabla$.  With the zeroth-order truncation
$f=g$, the Euler equations can be obtained. For the first-order
truncation
\begin{align*}
f=g-\tau (ug_x+vg_y+wg_z+g_t),
\end{align*}
the Navier-Stokes equations can be obtained.  More details on the
derivation of macroscopic equation can be found in
\cite{GKS-Xu1,GKS-Xu2}.

Originally, the high-order gas-kinetic scheme is developed in the
finite volume framework. In the following sections, the control
volume is hexahedron for simplicity. Integrating Eq.\eqref{hyper}
with respect to space, the semi-discretized finite volume scheme can
be expressed as
\begin{align}\label{semi}
\frac{\text{d} Q_i}{\text{d} t}=-\frac{1}{|\Omega_{i}|}\sum_{p=1}^6\iint_{\Sigma_{i_p}}\boldsymbol{F}(Q,t)\text{d}\sigma,
\end{align}
where $Q_i$ is the cell averaged conservative value of hexahedral
cell $\Omega_{i}$, $|\Omega_{i}|$ is the volume of $\Omega_{i}$,
$\Sigma_{i_p}$ is the quadrilateral common cell interface of
$\Omega_{i}$ and its neighboring cell $\Omega_{i_p}$. To achieve the
spatial accuracy, the Gaussian quadrature is used for the surface
integral of numerical fluxes over the cell interface
\begin{align*}
\iint_{\Sigma_{i_p}}\boldsymbol{F}(Q,t)\text{d}\sigma=\sum_{m_1,m_2}\omega_{m_1,m_2}F_{m_1,m_2}(t)V_{m_1,m_2},
\end{align*}
where $V_{m_1,m_2}$ is the area of cell interface related to
quadrature point. The numerical flux $F_{m_1,m_2}(t)$ at quadrature
point $\boldsymbol{x}_{m_1,m_2}$ in the global coordinate can be
given by
\begin{align*}
F_{m_1,m_2}(t)=\int\psi \boldsymbol{u}\cdot \boldsymbol{n}_{m_1,m_2}
f(\boldsymbol{x}_{m_1,m_2},t,\boldsymbol{u},\xi) \text{d}\Xi.
\end{align*}
In the computation, the numerical flux at quadrature points is
obtained in the local coordinate, and transferred to the global
coordinate. In the local coordinate, the gas distribution function
is constructed by the integral solution of BGK equation
Eq.\eqref{bgk-3d} as follows
\begin{equation*}
    f(\boldsymbol{x}_{m_1,m_2},t,\boldsymbol{u},\xi)=\frac{1}{\tau}\int_0^t
    g(\boldsymbol{x}',t',\boldsymbol{u},\xi)e^{-(t-t')/\tau}\text{d}t'+e^{-t/\tau}f_0(-\boldsymbol{u}t,\xi),
\end{equation*}
where $\boldsymbol{u}=(u,v,w)$ is the particle velocity in the local
coordinate,
$\boldsymbol{x}_{m_1,m_2}=\boldsymbol{x}'+\boldsymbol{u}(t-t')$ is
the trajectory of particles, $f_0$ is the initial gas distribution
function, and $g$ is the corresponding equilibrium state. With the
reconstruction of macroscopic variables, the second-order gas
distribution function at the cell interface can be expressed as
\begin{align}\label{flux}
    f(\boldsymbol{x}_{m_1,m_2},t,\boldsymbol{u},\xi)= & (1-e^{-t/\tau})g_0+((t+\tau)e^{-t/\tau}-\tau)(\overline{a}_1u+\overline{a}_2v+\overline{a}_3w)g_0\nonumber \\
    +                                                       & (t-\tau+\tau e^{-t/\tau}){\bar{A}} g_0\nonumber                                                            \\
    +                                                       & e^{-t/\tau}g_r[1-(\tau+t)(a_{1}^{r}u+a_{2}^{r}v+a_{3}^{r}w)-\tau A^r)](1-H(u))\nonumber                        \\
    +                                                       & e^{-t/\tau}g_l[1-(\tau+t)(a_{1}^{l}u+a_{2}^{l}v+a_{3}^{l}w)-\tau A^l)]H(u),
\end{align}
where $g_l, g_r$ are the equilibrium states  corresponding  to the
reconstructed variables $Q_l, Q_r$ at both sides of cell interface.
The coefficients in Eq.\eqref{flux} can be obtained by the
reconstructed directional derivatives and compatibility condition
\begin{align*}
    \langle a_{i}^{k}\rangle & =\frac{\partial Q_{k}}{\partial
        \boldsymbol{n}_i}, \langle
    a_{1}^{k}u+a_{2}^{k}v+a_{3}^{k}w+A^{k}\rangle=0,
\end{align*}
where $k=l,r$, $i=1,2,3$, the spatial derivatives
$\displaystyle\frac{\partial Q_{k}}{\partial { \boldsymbol{n}_i}}$
can be determined by spatial reconstruction and $\langle...\rangle$
are the moments of the equilibrium $g$ and defined by
\begin{align*}
    \langle...\rangle=\int g (...)\psi \text{d}\Xi.
\end{align*}
To avoid the extra reconstruction, the conservative variables
$Q_{0}$ and their spatial derivatives for equilibrium part can be
given by the compatibility condition Eq.\eqref{compatibility}
\begin{align*}
Q_0=&\int_{u>0}\psi g_{l}\text{d}\Xi+\int_{u<0}\psi g_{r}\text{d}\Xi,\\
\frac{\partial Q_{0}}{\partial {\boldsymbol{n}_i}}= & \int_{u>0}\psi
a_{i}^{l} g_{l}\text{d}\Xi+\int_{u<0}\psi a_{i}^{r}
g_{r}\text{d}\Xi.
\end{align*}
Similarly, the coefficients in Eq.\eqref{flux} for the equilibrium
part can be obtained by
\begin{align*}
\langle\overline{a}_i\rangle=\frac{\partial Q_{0}}{\partial
{\boldsymbol{n}_i}},
~\langle\overline{a}_1u+\overline{a}_2v+\overline{a}_3w+\overline{A}\rangle=0.
\end{align*}
Thus, the gas-distribution function at the cell interface has been
fully determined, and more details of the gas-kinetic scheme can be
found in \cite{GKS-Xu1}.

\section{Discontinuous Galerkin based gas-kinetic scheme}
In this paper, the three-dimensional structured meshes are considered
for simplicity. The computational domain $\Omega$ is divided into a
collection of non-overlap hexahedrons $ \mathcal{T}$, and each
cell is denoted as ${K_{ijk}=[x_{i-1/2}, x_{i+1/2}]\times[y_{j-1/2},
y_{j+1/2}]\times[z_{k+1/2}, z_{k+1/2}]}$, where $i=1,...,N_x$,
$j=1,...,N_y$ and $k=1,...,N_z$. To formulate
the discontinuous Galerkin method, the following Sobolev space
$\boldsymbol{V}^{k,m}$ is introduced
\begin{align*}
\boldsymbol{V}^{k,m}_h=\big\{\boldsymbol{v}_h\in [L^2(\Omega)]^m:
v_h|_{K_{ijk}}\in V^{k,m}\big\},
\end{align*}
where $v_h$ is one component for vector-valued $\boldsymbol{v}_h$,
$V^{k,m}$ is the space consisting of discontinuous polynomial
functions, and the degree equal or less than $k$ and $m=5$ is the
number of components. Multiplying the governing equation
Eq.\eqref{hyper}  by a test function $\boldsymbol{W}_h\in
\boldsymbol{V}^{k,m}_h$, integrating over $K_{ijk}$, and
performing an integration by parts, the following weak formulation
can be obtained
\begin{align}\label{DG}
    \frac{\text{d}}{\text{d} t}\int_{K_{ijk}} Q_h\cdot\boldsymbol{W}_h \text{d} \Omega
    +\sum_{\Gamma\in \partial K_{ijk}}\int_{\Gamma} \boldsymbol{F}(Q_h, t)\cdot \boldsymbol{n} \cdot \boldsymbol{W}_h\text{d}\Gamma
    -\int_{ K_{ijk}} \boldsymbol{F}(Q_h,t)\cdot \nabla\boldsymbol{W}_h\text{d}\Omega=0,
\end{align}
where $Q_h\in \boldsymbol{V}^{k,m}_h$ is the solution vector and
$\boldsymbol{n}$ is the unit outward normal vector of cell
interface. To discretize the weak formulation, two types of space
$\mathbb P_k$ and $\mathbb Q_k$ can be used for three-dimensional
flows, where $\mathbb{P}_k$ is the space of polynomials with the sum
of degree at most $k$ for the sum of all variables, and $\mathbb
Q_k$ is the space of tensor-product polynomials of degree at most
$k$ in each variable. For the three-dimensional computation, the
degrees of freedom is $(k+1)^3$ for $\mathbb Q_k$  and
$(k+1)(k+2)(k+3)/6$ for  $\mathbb P_k$. The space $\mathbb Q_k$
contains much more degrees of freedom than $\mathbb P_k$ and more
computational cost will be introduced, especially for
three-dimensional flows. To reduce the memory and improve
efficiency, $\mathbb P_k$ space is used. For the one-dimensional DG
scheme for the linear advection equation, the $(k+1)$-th order
optimal convergence and $(2k+1)$-th order super-convergence can be
obtained for the DG schemes with the $\mathbb P_k$ element. The
numerical results will show that the optimal convergence and
super-convergence can be achieved by DG-HGKS for three-dimensional
Euler equations as well.

The three-dimensional Legendre polynomials are chosen as the
basis functions of $\boldsymbol{V}^{k,m}_h$
\begin{align*}
\displaystyle \boldsymbol{B}_{\boldsymbol n}(\boldsymbol{x})=P_{i,n_x}(x)P_{j,n_y}(y)P_{k,n_z}(z),
\end{align*}
where $\boldsymbol n=(n_x, n_y, n_z)$ is the multi-index,
$P_{i,l}(x)=P_l(\xi)$  with $\xi=2(x-x_i)/\Delta x_i$ and $P_l(\xi)$
is the $l$th-order Legendre polynomial that are orthogonal on $[-1,
1]$. The time dependent solution vector $Q_h(\boldsymbol{x}, t)$ is
represented by the following expansion
\begin{align}\label{expansion}
Q_h(\boldsymbol{x},t)=\sum_{n=1}^NQ_{h, n}(t)\boldsymbol{B}_{\boldsymbol{n}}(\boldsymbol{x}),
\end{align}
where $Q_{h, n}(t)$ is the $n$-th unknown degree of freedom,
$n=n_x+n_y+n_z$ and $N$ is the total number of basis functions.
Taking
$\boldsymbol{W}_h=\boldsymbol{B}_{\boldsymbol{n}}(\boldsymbol{x})$
in Eq.\eqref{DG}, the DG formulation is equivalent to the following
system
\begin{align}\label{DG2}
\frac{\text{d}}{\text{d} t}\int_{K_{ijk}} Q_h\cdot
\boldsymbol{B}_{\boldsymbol n} \text{d} \Omega +\sum_{\Gamma\in
\partial K_{ijk}}\int_{\Gamma} \boldsymbol{F}_{b}(Q_h,t)\cdot
\boldsymbol{n} \cdot \boldsymbol{B}_{\boldsymbol n}\text{d}\Gamma
-\int_{K_{ijk}} \boldsymbol{F}_{v}(Q_h, t)\cdot
\nabla\boldsymbol{B}_{\boldsymbol n}\text{d}\Omega=0,
\end{align}
where $\boldsymbol{F}_{b}(Q_h,t)$ is the numerical flux across the
cell interface and $\boldsymbol{F}_{v}(Q_h,t)$ is the numerical flux
inside each volume.  According to the expansion
Eq.\eqref{expansion}, the DG formulation Eq.\eqref{DG2} leads to a
system of ordinary differential equations
\begin{align*}
\frac{\text{d}Q_{h}}{\text{d}t}=\boldsymbol{M}^{-1}\mathcal{R}(Q_h)=\mathcal {L}(Q_h),
\end{align*}
where each term of the mass matrix $\boldsymbol{M}_{N\times N}$ is
defined by
\begin{align*}
\boldsymbol{M}_{n_1,n_2}= \int_{K_{ijk}} \boldsymbol{B}_{\boldsymbol{n}_1}
\boldsymbol{B}_{\boldsymbol {n}_2}\text{d} \Omega,
\end{align*}
and the operator $\mathcal{R}(Q_h)$ is given by
\begin{align}\label{operator}
\mathcal{R}(Q_h)=-\sum_{\Gamma\in \partial
K_{ijk}}\int_{\Gamma} \boldsymbol{F}_{b}(Q_h, t)\cdot
\boldsymbol{n} \cdot \boldsymbol{B}_{\boldsymbol n}\text{d}S
+\int_{K_{ijk}} \boldsymbol{F}_v(Q_h, t)\cdot
\nabla\boldsymbol{B}_{\boldsymbol n}\text{d}\Omega.
\end{align}
The Gaussian quadrature is used for the surface integral
\begin{align}\label{surface}
\int_{\Gamma} \boldsymbol{F}_{b}(\boldsymbol{U}_h,t)\cdot
\boldsymbol{n}  \cdot \boldsymbol{B}_{\boldsymbol{n}}\text{d}S
=\sum_{m_1,m_2}\omega_{m_1,m_2}
\boldsymbol{F}_{b}(\boldsymbol{x}_{m_1,m_2},t)\cdot
\boldsymbol{n} \cdot
\boldsymbol{B}_{\boldsymbol{n}}(\boldsymbol{x}_{m_1,m_2})|\Gamma|,
\end{align}
where $\omega_{m_1,m_2}$ is the quadrature weight for 2D Gaussian
quadrature, $|\Gamma|$ is the area of $\Gamma$, and the numerical
flux $\boldsymbol{F}_{b}(\boldsymbol{x}_{m_1,m_2},t)$ can be
obtained by taking moments of gas distribution function given by
Eq.\eqref{flux}
\begin{align*}
\boldsymbol{F}_{b}(\boldsymbol{x}_{m_1,m_2},t)=\int\psi
\boldsymbol{u}\cdot \boldsymbol{n}_{m_1,m_2}
f(\boldsymbol{x}_{m_1,m_2},t,\boldsymbol{u},\xi)\text{d}\Xi.
\end{align*}
Similarly, the Gaussian quadrature is also used for the  volume
integral
\begin{align}\label{volume}
\int_{K_{ijk}}\boldsymbol{F}_v(Q_h) \cdot
\nabla\boldsymbol{B}_i\text{d}\Omega=
\sum_{m_1,m_2,m_3}\omega_{m_1,m_2,m_3}\boldsymbol{F}_v(\boldsymbol{x}_{m_1,m_2,m_3},t)\cdot
\nabla\boldsymbol{B}_i(\boldsymbol{x}_{m_1,m_2,m_3})|K_{ijk}|,
\end{align}
where $\omega_{m_1,m_2,m_3}$ is the quadrature weight for 3D
Gaussian quadrature, $|K_{ijk}|$ is the volume of $K_{ijk}$, and the
numerical flux $\boldsymbol{F}_v(\boldsymbol{x}_{m_1,m_2,m_3},t)$
can be given by
\begin{align*}
\boldsymbol{F}_{v}(\boldsymbol{x}_{m_1,m_2,m_3},t)=\int\psi
\boldsymbol{u}\cdot \boldsymbol{n}_{m_1,m_2,m_3}
f(\boldsymbol{x}_{m_1,m_2,m_3},t,\boldsymbol{u},\xi)\text{d}\Xi.
\end{align*}

\begin{rmk}
The macroscopic flow variables inside a cell are updated
automatically with smooth assumption, and the simplified version of
gas distribution function is used
\begin{align}\label{flux-2}
f(\boldsymbol{x}_{m_1,m_2,m_3},t,\boldsymbol{u},\xi)=g_0\big(1-\tau(a_1u+a_2v+a_3w+
A) +At\big),
\end{align}
where the coefficients in Eq.\eqref{flux-2} can be obtained by the
reconstructed directional derivatives and compatibility condition
\begin{align*}
\langle a_{i}\rangle & =\frac{\partial Q}{\partial\boldsymbol{n}_i},
\langle a_{1}u+a_{2}v+a_{3}w+A\rangle=0.
\end{align*}
In the numerical cases, the smooth flow without discontinuities are
considered and the collision time takes
\begin{align*}
\tau=\frac{\mu}{p}.
\end{align*}
For the accuracy tests, the inviscid flows are considered. The
collision time $\tau$ and viscous coefficient $\mu$ take $0$, and
the gas distribution functions for the surface and volume
integrations reduce to
\begin{align*}
f=g_0\big(1+At\big).
\end{align*}
In this paper, we mainly focus on the smooth flows. For the flows
with discontinuities, and the limiters and trouble indicators need
to be used, the coupling with the current scheme will be
investigated in the future.
\end{rmk}
\begin{rmk}
To achieve the spatial accuracy, $\mathbb P_2$ and $\mathbb P_3$
spaces are used. The two-point and three-point Gaussian quadratures
are adopted, i.e. For the $\mathbb P_2$ space, $2^2$ numerical
fluxes are calculated for surface integration Eq.\eqref{surface} and
$2^3$ numerical fluxes are calculated for volume integration
Eq.\eqref{volume}; For the $\mathbb P_3$ space, $3^2$ and $3^3$
numerical fluxes are calculated for the surface and volume
integrations respectively.
\end{rmk}

In the classical numerical methods, the numerical fluxes are usually
provided by the approximate or exact Riemann solvers
\cite{Riemann-appro},  and the Runge-Kutta method is used for
temporal accuracy \cite{TVD-RK}.  A two-stage fourth-order
time-accurate discretization was developed for the hyperbolic
equations with the generalized Riemann problem (GRP) solver
\cite{GRP-high-1, GRP-high-2} and the gas-kinetic scheme
\cite{GKS-high-1}. Consider the following time-dependent equation
\begin{align*}
\frac{\text{d}Q}{\text{d} t}=\mathcal {L}(Q),
\end{align*}
with the initial condition at $t_n$,
\begin{align*}
Q(t)|_{t=t_n}=Q^n,
\end{align*}
where $\mathcal {L}(Q)=\boldsymbol{M}^{-1}\mathcal{R}(Q)$.
Introducing an intermediate state at $t^*=t_n+\Delta t/2$,
\begin{equation}\label{two-stage}
\begin{split}
&Q^*=Q^n+\frac{1}{2}\Delta t\mathcal{L}(Q^n)+\frac{1}{8}\Delta t^2\mathcal{L}_t(Q^n),\\
Q^{n+1}&=Q^n+\Delta t\mathcal {L}(Q^n)+\frac{1}{6}\Delta
t^2\big(\mathcal{L}_t(Q^n)+2\mathcal{L}_t(Q^*)\big).
\end{split}
\end{equation}
It can be proved that for hyperbolic equations the two-stage time
stepping method Eq.\eqref{two-stage} provides a fourth-order time
accurate solution for $Q(t)$ at $t=t_n+\Delta t$. To implement the
two-stage fourth-order method, the following terms need to be
constructed by the Gaussian quadrature point
\begin{align*}
\boldsymbol{F}_{b}(Q_h(\boldsymbol{x}_{m_1,m_2},t_n)), &~\partial_t\boldsymbol{F}_{b}(Q_h(\boldsymbol{x}_{m_1,m_2},t_n)),\\
\boldsymbol{F}_{v}(Q_h(\boldsymbol{x}_{m_1,m_2,m_3},t_n)),
&~\partial_t\boldsymbol{F}_{v}(Q_h(\boldsymbol{x}_{m_1,m_2,m_3},t_n)).
\end{align*}
A linear function is used to approximate the time dependent
numerical flux
\begin{align}\label{expansion-1}
\boldsymbol{F}(Q_h(\boldsymbol{x},t))=\boldsymbol{F}(Q_h(\boldsymbol{x}, t_n))+\partial_t\boldsymbol{F}(Q_h(\boldsymbol{x},t_n))(t-t_n),
\end{align}
where the subscripts $b$, $v$ and indexes for quadrature points are
omitted. Integrating Eq.\eqref{expansion-1} over $[t_n, t_n+\Delta
t]$ and $[t_n, t_n+\Delta t/2]$, we have the following two equations
\begin{align*}
\boldsymbol{F}(Q_h(\boldsymbol{x}, t_n))\Delta t&
+\frac{1}{2}\partial_t\boldsymbol{F}(Q_h(\boldsymbol{x}, t_n))\Delta t^2 =\int_{t_n}^{t_n+\delta}\boldsymbol{F}(Q_h(\boldsymbol{x},t))\text{d}t,\\
\frac{1}{2}\boldsymbol{F}(Q_h(\boldsymbol{x}, t_n))\Delta t &
+\frac{1}{8}\partial_t\boldsymbol{F}(Q_h(\boldsymbol{x}, t_n))\Delta t^2 =\int_{t_n}^{t_n+\delta/2}\boldsymbol{F}(Q_h(\boldsymbol{x},t))\text{d}t.
\end{align*}
The coefficients can be determined by solving the linear system, and
$\mathcal {L}(Q^n)$ and $\mathcal{L}_t(Q^n)$ can be given.
Similarly, $\mathcal {L}(Q^*)$ and $\mathcal{L}_t(Q^*)$ can be given
at $t=t_n+\Delta/2$. Thus, the DG-HGKS can be implemented.

\section{GPU architecture and code design}
Graphics processing unit (GPU)  is a form of hardware acceleration,
which is originally developed for graphics manipulation and is
extremely efficient at processing large amounts of data in parallel,
such as  pixel transformations. Currently, GPU has gained
significant popularity in large-scale scientific computations \cite{GPU}. Due
to the different purposes of designed, the architectures of CPU and
GPU are different. A central processing unit (CPU) is designed to
execute a sequence of tasks quickly, but is limited in the number of
threads which are handled in parallel. Meanwhile, a GPU is designed
to execute highly-parallel computing tasks, and it may handle
thousands of threads at the same time.

\begin{figure}[!h]
    \centering
    \includegraphics[width=0.65\textwidth]{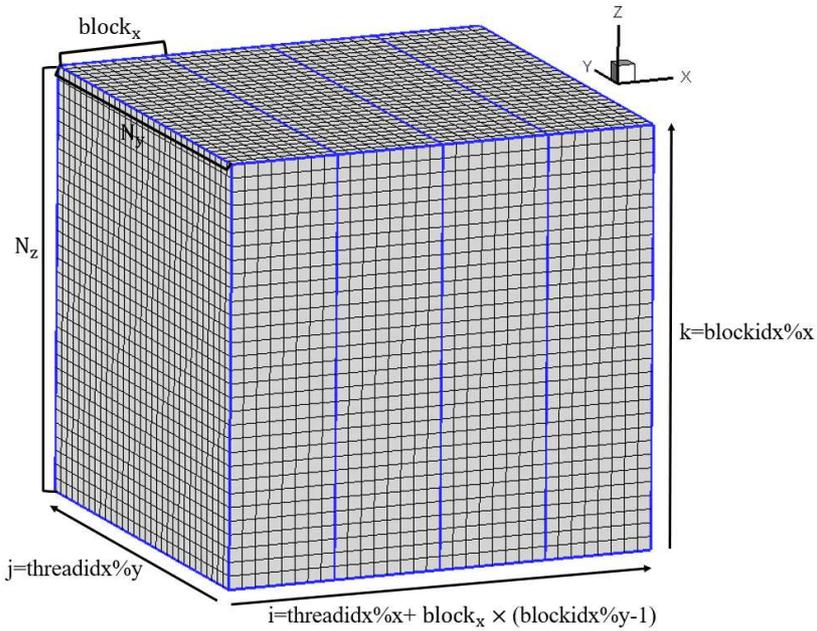}
    \caption{\label{thread-cell} Connection between threads and cells.}
\end{figure}

\begin{table}[!h]
    \begin{center}
    \def\temptablewidth{0.9\textwidth} {\rule{\temptablewidth}{1.0pt}}
    \begin{tabular*}{\temptablewidth}{@{\extracolsep{\fill}}c|c|c|c|c}
    ~~              & threadidx\%x  & threadidx\%y & blockidx\%x & blockidx\%y  \\
    \hline
    range of value & $[1, \text{block}_x]$ & $[1, N_y]$ & $[1, N_z]$ & $[1,  N_x/\text{block}_x]$ \\
    \end{tabular*}
    {\rule{\temptablewidth}{1.0pt}}
    \end{center}
    \caption{\label{tab-variable-range} Range of values for "threadidx" and "blockidx".}
    \end{table}

In this paper, to accelerate the computation, DG-HKGS is implemented
on GPU using CUDA. The way to
implement the algorithm in parallel on GPU is using kernels. The
kernel is a subroutine, which executes at the same time by many
threads on GPU. According to the algorithm of DG-HGKS proposed
above, the subroutines of surface integrals for $\boldsymbol{F}_{v}$
and $\partial_t\boldsymbol{F}_{v}$ in Eq.\eqref{surface}, the volume
integrals for $\boldsymbol{F}_{b}$ and
$\partial_t\boldsymbol{F}_{b}$ in Eq.\eqref{volume} and the update
of flow variables $Q_{h}$ can be implemented parallelly. Therefore,
these subroutines can be coded as kernels. In order to improve the
efficiency, $\boldsymbol{F}_{v}$, $\partial_t\boldsymbol{F}_{v}$,
$\boldsymbol{F}_{b}$, $\partial_t\boldsymbol{F}_{b}$ and $Q_{h}$ are
saved as global variables stored in global memory. The code can be
expressed with following steps
\begin{enumerate}
\item dispatching $Q_{h}(t_n)$ of $K_{ijk}$ and its neighboring cells from global memory to local at $t_n$,
\item computing $\boldsymbol{F}_{v}$, $\partial_t\boldsymbol{F}_{v}$ and $\boldsymbol{F}_{b}$, $\partial_t\boldsymbol{F}_{b}$,
and obtaining $\mathcal{R}(Q_h)$ by  Eq.\eqref{operator} at  $t_n$,
\item updating $Q_{h}(t^*)$ of cell $K_{ijk}$ according to Eq.\eqref{two-stage},
\item dispatching the updated $Q_{h}(t^*)$ back to global memory,
\item repeating the procedures above at the intermediate stage $t^*$.
\end{enumerate}
The CUDA threads are organized into thread blocks, and thread blocks
constitute a grid which may be seen as a computational structure on
GPUs. In this study, the three-dimensional structured meshes are
used, a simple partition of grid is established, which is shown in
Figure.\ref{thread-cell}. Assume that the total number of meshes is
$N_x\times N_y\times N_z$, and the three-dimensional  computational
domain is divided into $\text{block}_x$ parts in $x$-direction,
where $\text{block}_x$ is an integer defined according to tests. The
variables "dimGrid" and "dimBlock" are defined to set a
two-dimensional grid, which consist of  two-dimensional blocks as
follows
\begin{align*}
     &  {\rm dimGrid}={\rm dim3}(N_z, N_x/\text{block}_x, 1), \\
     &  {\rm dimBlock}={\rm dim3}(\text{block}_x,N_y,1).
\end{align*}
Each thread gets a block index "blockidx" and a thread index
"threadidx", where the two variables have two components and the
range of values are given in Table.\ref{tab-variable-range}. It is
natural to assign one thread to complete computations of a cell
$K_{ijk}$. The one-to-one correspondence of thread block index and
cell index $(i,j,k)$ is given as follows
\begin{align*}
i&=\text{threadidx}\%\text{x}+\text{block}_x*(\text{blockidx}\%\text{y}-1),\\
j&=\text{threadidx}\%\text{y}, \\
k&= \text{blockidx}\%\text{x}.
\end{align*}
The Nvidia GPU is consisted of multiple streaming multiprocessors
(SMs). Each block of grid is distributed to one SM, and the threads
of block are executed in parallel on SM. The executions are
implemented automatically by GPU. Thus, the GPU code can be
implemented with specifying kernels and grids, and the performance
will be tested.

\section{Numerical simulation and discussion}
\subsection{Accuracy tests}
For the one-dimensional scalar hyperbolic equation and hyperbolic
system, the optimal $(k+1)$-th order convergence can be obtained
with $\mathbb P_k$ elements. For the multidimensional cases, the
optimal $(k+1)$-th order can be achieved, when the piecewise tensor
product polynomials are used, i.e., $\mathbb Q_k$ elements.
Recently,  for the two-dimensional scaler hyperbolic equations, the
optimal $(k+1)$-th order convergence is proved for $\mathbb P_k$
elements on uniform Cartesian meshes \cite{DG-convergence}. For
one-dimensional linear hyperbolic equations, with  suitable
correction function, the super-convergence
\cite{DG-superconvergence}  can be proved  at the downwind points
and for the domain average  with quasi-uniform meshes and some
suitable initial discretization.  This technique was developed to
study other kinds of DG methods, such as the local DG method, the
direct DG method. Two cases will be provided to test the optimal
convergence and super-convergence of DG-HGKS for nonlinear Euler
equations. To test the optimal convergence, the $L^1$ and $L^2$
errors are calculated, and the definitions are given as follows
\begin{align*}
e_{L^1}&=\Big( \sum_{\Omega_{e} \in \mathcal{T}_{h}}\big( \int_{\Omega_{e}}|u-u_{h}|d\Omega\big)\Big),\\
e_{L^2}&=\Big( \sum_{\Omega_{e} \in \mathcal{T}_{h}}\big( \int_{\Omega_{e}}|u-u_{h}|^{2}d\Omega\big)\Big)^{1/2},
\end{align*}
where $u$ is the exact solution and $u_{h}$ is the approximate
numerical solution, and Gaussian quadrature is used for the
integration. To test super-convergence of DG scheme, the error for
the cell averaged value is defined as follows
\begin{align*}
e_{c}=\Big( \sum_{\Omega_{e} \in \mathcal{T}_{h}}\big( \int_{\Omega_{e}}|\bar{u}-\bar{u}_{h}|^{2}d\Omega\big)\Big)^{1/2}.
\end{align*}
where $\bar{u}$ is cell averaged value of exact solution and
$\bar{u}_{h}$ is cell averaged value of numerical solution.

The first case is the advection of density perturbation for the
two-dimensional and three-dimensional flows, which are linear
problems for accuracy test. For the two-dimensional case, the
physical domain is $[0,2]\times[0,2]$ and the initial condition is
set as follows
\begin{align*}
\rho_0(x, y)=1+0.2\sin(\pi(x+y)),~p_0(x,y)=1,~U_0(x,y)=1,~V_0(x,y)=1.
\end{align*}
The periodic boundary conditions are applied at boundaries, and the
exact solution is
\begin{align*}
\rho(x,y,t)=1+0.2\sin(\pi(x+y-t)),~p(x,y,t)=1,~U(x,y,t)=1,~V(x,y,t)=1.
\end{align*}
This case is tested with uniform and nonuniform meshes. For the
uniform meshes,  the mesh with $\Delta x=\Delta y=2/N$ are tested.
For the nonuniform mesh, the gird points are given by the following
transformation
\begin{align*}
\begin{cases}
\displaystyle x=\xi+0.05\sin (\pi \xi),\\
\displaystyle y=\eta+0.05\sin (\pi \eta).
\end{cases}
\end{align*}
For the three-dimensional case, the physical domain is
$[0,2]\times[0,2]\times[0,2]$ and the initial condition is set as
follows
\begin{align*}
\rho_0 & (x, y, z)=1+0.2\sin(\pi(x+y+z)),~p_0(x,y,z)=1, \\
           & U_0(x,y,z)=1,~V_0(x,y,z)=1,~W_0(x,y,z)=1.
\end{align*}
The periodic boundary conditions are applied at boundaries, and the
exact solution is
\begin{align*}
    \rho(x,y & ,z,t)=1+0.2\sin(\pi(x+y+z-t)),~p(x,y,z,t)=1, \\
             & U(x,y,z,t)=1,~V(x,y,z,t)=1,~W(x,y,z,t)=1.
\end{align*}
This case is also tested with uniform and nonuniform meshes. For the
uniform meshes, the mesh with $\Delta x=\Delta y=\Delta z=2/N$  are
tested. For the nonuniform mesh, the gird point is given by the
following transformation
\begin{align*}
\begin{cases}
\displaystyle x=\xi+0.05\sin (\pi \xi),\\
\displaystyle y=\eta+0.05\sin (\pi \eta),\\
\displaystyle z=\zeta+0.05\sin (\pi \zeta).
\end{cases}
\end{align*}
In the computation, the $\mathbb P_2$ and $\mathbb P_3$ elements are
tested. The $L^1$ and $L^2$ errors and orders of accuracy for the
two-dimensional problem are given in Table.\ref{tab-2d-1},
Table.\ref{tab-2d-2}, Table.\ref{tab-2d-3} and Table.\ref{tab-2d-4}
for both uniform and nonuniform meshes, and the expected $(k+1)$-th
order of accuracy is observed. For the hyperbolic scalar equation,
$(2k+1)$-th order super-convergence is expected
\cite{DG-superconvergence} for the error of cell averaged value. For
the current method, the error of cell averaged value and orders of
accuracy are also given in Table.\ref{tab-2d-1},
Table.\ref{tab-2d-2}, Table.\ref{tab-2d-3} and Table.\ref{tab-2d-4}
for both uniform and nonuniform meshes, and the $2k$-th order
super-convergence is observed for $\mathbb P_k$ elements. It might
be because the fewer Gaussian quadrature points for
Eq.\eqref{surface} and Eq.\eqref{volume}. The $L^1$ and $L^2$ errors
and orders of accuracy for the three-dimensional problem are given
in Table.\ref{tab-3d-1}, Table.\ref{tab-3d-2}, Table.\ref{tab-3d-3}
and Table.\ref{tab-3d-4}, and the expected $(k+1)$-th order of
accuracy are observed as well. Due to the very small time step, the
super-convergence for three-dimensional problem is not given. It is
also believed that the $2k$-th order super-convergence can be
achieved for $\mathbb P_k$ elements.

\begin{table}[!h]
\begin{center}
\def\temptablewidth{0.95\textwidth} {\rule{\temptablewidth}{1.0pt}}
\begin{tabular*}{\temptablewidth}{@{\extracolsep{\fill}}c|cc|cc|cc}
mesh     & $L^1$ error  &    Order    &  $L^2$ error &  Order   & $e_{u, c}$ error  &    Order \\
\hline
$8^2$          &   1.2632E-02 & ~~      & 5.0047E-03 &   ~~         &   5.6223E-04 & ~~        \\
$16^2$         &   1.2982E-03 & 3.2824  & 5.1162E-04 &  3.2902      &   2.4338E-05 & 4.5299    \\
$32^2$         &   1.5215E-04 & 3.0930  & 5.9793E-05 &  3.0970      &   1.1353E-06 & 4.4220    \\
$64^2$         &   1.8633E-05 & 3.0296  & 7.3140E-06 &  3.0312      &   6.1540E-08 & 4.2054    \\
\end{tabular*}
{\rule{\temptablewidth}{1.0pt}}
\end{center}
\caption{\label{tab-2d-1} Accuracy test: 2D advection of density
perturbation for $\mathbb P_2$ element with uniform cells.}
\begin{center}
\def\temptablewidth{0.95\textwidth}
{\rule{\temptablewidth}{1.0pt}}
\begin{tabular*}{\temptablewidth}{@{\extracolsep{\fill}}c|cc|cc|cc}
mesh     &   $L^1$ error  &Order &   $L^2$ error  &Order  & $e_{u, c}$ error  &    Order  \\
\hline
$8^2$          &   6.4334E-04 & ~~      & 3.1235E-04 &   ~~     &   4.5650E-06 & ~~       \\
$16^2$         &   3.8232E-05 & 4.0727  & 1.7848E-05 &  4.1293  &   4.5357E-08 & 6.6531   \\
$32^2$         &   2.3296E-06 & 4.0366  & 1.0881E-06 &  4.0359  &   7.0660E-10 & 6.0043   \\
$64^2$         &   1.4433E-07 & 4.0126  & 6.7525E-08 &  4.0102  &   1.0552E-11 & 6.0653   \\
\end{tabular*}
{\rule{\temptablewidth}{1.0pt}}
\end{center}
\caption{\label{tab-2d-2} Accuracy test: 2D advection of density
perturbation for $\mathbb P_3$ element with uniform cells.}
\begin{center}
\def\temptablewidth{0.95\textwidth}{\rule{\temptablewidth}{1.0pt}}
\begin{tabular*}{\temptablewidth}{@{\extracolsep{\fill}}c|cc|cc|cc}
mesh     &   $L^1$ error  &Order &   $L^2$ error  &Order   & $e_{u, c}$ error  &    Order     \\
\hline
$8^2$          &   1.3484E-02  & ~~      & 5.4209E-03  &   ~~      & 6.2038E-04  &   ~~       \\
$16^2$         &   1.3756E-03  & 3.2932  & 5.4617E-04  &   3.3111  & 2.7181E-05  &   4.5125   \\
$32^2$         &   1.6055E-04  & 3.0989  & 6.3559E-05  &   3.1032  & 1.2526E-06  &   4.4396   \\
$64^2$         &   1.9643E-05  & 3.0310  & 7.7662E-06  &   3.0328  & 6.6467e-08  &   4.2361   \\
\end{tabular*}
{\rule{\temptablewidth}{1.0pt}}
\end{center}
\caption{\label{tab-2d-3} Accuracy test: 2D advection of density
perturbation for $\mathbb P_2$ element with nonuniform cells.}
\begin{center}
\def\temptablewidth{0.95\textwidth}
{\rule{\temptablewidth}{1.0pt}}
\begin{tabular*}{\temptablewidth}{@{\extracolsep{\fill}}c|cc|cc|cc}
mesh     &   $L^1$ error  &Order &   $L^2$ error  &Order  & $e_{u, c}$ error  &    Order  \\
\hline
$8^2$          &   8.3231E-04  & ~~      & 3.6625E-04  &   ~~      & 6.0414E-06  &   ~~    \\
$16^2$         &   4.4123E-05  & 4.2375  & 1.9564E-05  &   4.2266  & 8.4484E-08  &   6.1601   \\
$32^2$         &   2.6616E-06  & 4.0512  & 1.1505E-06  &   4.0878  & 5.8786E-10  &   7.1671   \\
$64^2$         &   1.6712E-07  & 3.9933  & 7.0616E-08  &   4.0262  & 1.1928E-11  &   5.6230   \\
\end{tabular*}
{\rule{\temptablewidth}{1.0pt}}
\end{center}
\caption{\label{tab-2d-4} Accuracy test: 2D advection of density
perturbation for $\mathbb P_3$ element with nonuniform cells.}
\end{table}

\begin{table}[!h]
\begin{center}
\def\temptablewidth{0.75\textwidth} {\rule{\temptablewidth}{1.0pt}}
\begin{tabular*}{\temptablewidth}{@{\extracolsep{\fill}}c|cc|cc}
mesh     & $L^1$ error  &    Order    &  $L^2$ error &  Order     \\
\hline
$8^3$          &   3.6718E-02 & ~~      & 1.4834E-03 &   ~~       \\
$16^3$         &   2.9226E-03 & 3.6512  & 1.1441E-03 &  3.6965    \\
$32^3$         &   2.9518E-04 & 3.3076  & 1.1712E-04 &  3.2882    \\
$64^3$         &   3.4117E-05 & 3.1130  & 1.3830E-05 &  3.0821    \\
\end{tabular*}
{\rule{\temptablewidth}{1.0pt}}
\end{center}
\caption{\label{tab-3d-1} Accuracy test:  3D advection of density
perturbation for $\mathbb P_2$ element with uniform cells.}
\begin{center}
\def\temptablewidth{0.75\textwidth}
{\rule{\temptablewidth}{1.0pt}}
\begin{tabular*}{\temptablewidth}{@{\extracolsep{\fill}}c|cc|cc}
mesh     &   $L^1$ error  &Order &   $L^2$ error  &Order    \\
\hline
$8^3$          &   1.6627E-03 & ~~      & 8.9026E-04 &   ~~       \\
$16^3$         &   9.0567E-05 & 4.1984  & 5.3092E-05 &   4.0677   \\
$32^3$         &   5.6787E-06 & 3.9954  & 3.3471E-06 &  3.9875    \\
$64^3$         &   3.5448E-07 & 4.0018  & 2.0974E-07 &  3.9962    \\
\end{tabular*}
{\rule{\temptablewidth}{1.0pt}}
\end{center}
\caption{\label{tab-3d-2} Accuracy test:  3D advection of density
perturbation for $\mathbb P_3$ element with uniform cells.}
\begin{center}
\def\temptablewidth{0.75\textwidth}{\rule{\temptablewidth}{1.0pt}}
\begin{tabular*}{\temptablewidth}{@{\extracolsep{\fill}}c|cc|cc}
mesh     &   $L^1$ error  &Order &   $L^2$ error  &Order       \\
\hline
$8^3$          &   3.9464E-02  & ~~          & 1.6035E-02  &   ~~          \\
$16^3$         &   3.1274E-03  & 3.6575      & 1.2335E-03  &   3.7004      \\
$32^3$         &   3.1367E-04  & 3.3177      & 1.2516E-04  &   3.3009      \\
 $64^3$        &   3.6132E-05  & 3.1179      & 1.4737E-05  &   3.0864      \\
\end{tabular*}
{\rule{\temptablewidth}{1.0pt}}
\end{center}
\caption{\label{tab-3d-3} Accuracy test:  3D advection of density
perturbation for $\mathbb P_2$ element with nonuniform cells.}
\begin{center}
\def\temptablewidth{0.75\textwidth}
{\rule{\temptablewidth}{1.0pt}}
\begin{tabular*}{\temptablewidth}{@{\extracolsep{\fill}}c|cc|cc}
mesh     &   $L^1$ error  &Order &   $L^2$ error  &Order   \\
\hline
$8^3$         &   2.1888E-03  & ~~      & 1.0100E-03  &   ~~         \\
$16^3$        &   1.0973E-04  & 4.3181  & 5.4167E-05  &   4.2208     \\
$32^3$        &   6.7226E-06  & 4.0288  & 3.3516E-06  &   4.0145     \\
 $64^3$       &   4.1725E-07  & 4.0100  & 2.0926E-07  &   4.0015     \\
\end{tabular*}
 {\rule{\temptablewidth}{1.0pt}}
\end{center}
\caption{\label{tab-3d-4} Accuracy test:  3D advection of density
perturbation for $\mathbb P_3$ element with nonuniform cells.}
\end{table}

\begin{table}[!h]
\begin{center}
\def\temptablewidth{0.95\textwidth}{\rule{\temptablewidth}{1.0pt}}
\begin{tabular*}{\temptablewidth}{@{\extracolsep{\fill}}c|cc|cc|cc}
mesh     & $L^1$ error  &    Order    &  $L^2$ error &  Order   & $e_{c}$ error  &    Order   \\
\hline
$20^2$          &   9.9414E-01   & ~~      & 4.7547E-02   &   ~~      & 2.5500E-04  &   ~~   \\
$40^2$          &   1.2966E-01   & 2.9387  & 7.0093E-03   &  2.7620   & 1.5981E-05  &   3.9961   \\
$80^2$          &   1.8311E-02   & 2.8240  & 1.1796E-03   &  2.5709   & 1.4020E-06  &   3.5107   \\
$160^2$         &   2.3769E-03   & 2.9456  & 2.0082E-04   &  2.5544   & 1.3725E-07  &   3.3527   \\
\end{tabular*}
{\rule{\temptablewidth}{1.0pt}}
\end{center}
\caption{\label{tab-iso-vor-1} Accuracy test: isotropic vortex
propagation problem for $\mathbb P_2$ element with uniform cells.}
\begin{center}
\def\temptablewidth{0.95\textwidth}
{\rule{\temptablewidth}{1.0pt}}
\begin{tabular*}{\temptablewidth}{@{\extracolsep{\fill}}c|cc|cc|cc}
mesh     & $L^1$ error  &    Order    &  $L^2$ error &  Order    & $e_{c}$ error  &    Order  \\
\hline
$20^2$         &   8.7222E-02 & ~~         & 4.0368E-03 &   ~~   & 1.2223E-05  &   ~~  \\
$40^2$         &   3.6128E-03 & 4.5935  & 2.0043E-04 &  4.3320   & 3.0863E-07  &   5.3076   \\
$80^2$         &   1.9642E-04 & 4.2011  & 1.2460E-05 &  4.0077   & 1.0386E-08  &   4.8931   \\
$160^2$        &   1.1337E-05 & 4.1148  & 8.0350E-07 &  3.9549   & 4.1164E-10  &   4.6572   \\
\end{tabular*}
{\rule{\temptablewidth}{1.0pt}}
\end{center}
\caption{\label{tab-iso-vor-2} Accuracy test: isotropic vortex
propagation problem for $\mathbb P_3$ element with uniform cells.}
\end{table}

The second case is the isotropic vortex propagation problem, and it
is a non-linear case for accuracy test. The initial condition is
given by a mean flow $(\rho, U, V,  p) = (1, 1, 1, 1)$. An isotropic
vortex is added to the mean flow, i.e., with perturbation in $U$,
$V$ and temperature $T = p/\rho$, and no perturbation in entropy $S
=p/\rho^\gamma$. The perturbation is given by
\begin{align*}
(\delta U,\delta V)=\frac{\epsilon}{2\pi}e^{\frac{(1-r^2)}{2}}(-y,x), ~ \delta
T=-\frac{(\gamma-1)\epsilon^2}{8\gamma\pi^2}e^{1-r^2},   ~\delta S=0,
\end{align*}
where $r^2=x^2+y^2$ and the vortex strength $\epsilon=5$. The
computational domain is $[0, 10] \times [0, 10]$, and the periodic
boundary conditions are imposed on the boundaries in all directions.
The exact solution is the perturbation of isotropic vortex which
propagates with the velocity $(1,1)$. In the computation, this case
is tested by $\mathbb P_2$ and $\mathbb P_3$ elements, and the
uniform mesh with $N^2$ cells are used. The $L^1$ and $L^2$ errors
and orders of accuracy are given in Table.\ref{tab-iso-vor-1} and
Table.\ref{tab-iso-vor-2}, and the expected $(k+1)$-th order of
accuracy is observed for $\mathbb P_k$ elements. The errors for cell
averaged value are given in Table.\ref{tab-iso-vor-1} and
Table.\ref{tab-iso-vor-2}, and the $(k+3/2)$-th order
super-convergence is observed.

In these two cases, the optimal order of accuracy are obtained for
both linear and nonlinear problems for $\mathbb P_k$ elements for
two and three-dimensional nonlinear hyperbolic system. Meanwhile,
the super-convergence are also observed for both linear and
nonlinear problems. For the nonlinear systems, the error estimates
of optimal convergence and super-convergence for $\mathbb P_k$
elements are still very challenging for DG schemes, and more studies
will be investigated in the future.

\subsection{Efficiency comparison of CPU and GPU}
The efficiency comparison of CPU and GPU code is provided with the
three-dimensional advection of density perturbation. The CPU code is
computed with  an octa-core  Intel i7-9700 CPU using Intel Fortran
compiler with OpenMP directives, while Nvidia TITAN RTX and Nvidia
Tesla V100 are used for GPU computation with Nvidia CUDA and
NVFORTRAN compiler. The detailed parameters of GPU and CPU are given
in Table.\ref{GPU-CPU-1}. Double precision is used in computation.
The number of the uniform meshes used for tests is from $16^3$ to
$128^3$, and the computational times are recorded at $t=2$.  The
computational times and speedups for the DG-HGKS with different
meshes are given in Table.\ref{GPU-CPU-2}, where speedup is defined
as
\begin{equation*}
\displaystyle \text{speedup}=\frac{\text{computational time of CPU}}{\text{computational time of GPU}}.
\end{equation*}
Theoretically, the total computation amount increases by a factor of
16, when the number of cells doubles in every direction.  As
expected, the computational time for the parallel CPU code increase
by the same factor approximately. The GPU computational time
increases by a factor less than 16  for coarse meshes, and the
factor reaches 16 for the finer meshes. The speedup increases when
the number of cells increases, and reaches an almost constant value.
Compared with the CPU code,  6x-7x speedup is achieved for TITAN RTX
and 10x-11x speedup is achieved for Tesla V100. In the future, the
code of HGKS will be developed with multiple GPUs, and  more
challenging problems for compressible flows will be investigated.

\begin{table}[!h]
\begin{center}
\def\temptablewidth{0.95\textwidth}{\rule{\temptablewidth}{1.0pt}}
\begin{tabular*}{\temptablewidth}{@{\extracolsep{\fill}}c|c|c|c}
~  &  Intel i7-9700 CPU    & Nvidia TITAN RTX & Nvidia Tesla V100  \\
\hline
clock rate              & 3.0 GHz      &1.77 GHz        &1.53 GHz \\
\hline
stream multiprocessor   &-             &72              &80  \\
\hline
FP64 core per SM        &-             &2               &32 \\
\hline
double precision        &384 GFLOPS    &509.8 GFLOPS    &7.834 TFLOPS \\
\end{tabular*}
{\rule{\temptablewidth}{1.0pt}}
\end{center}
\caption{\label{GPU-CPU-1} Efficiency comparison: the detailed
parameters of GPU and CPU.}
\begin{center}
\def\temptablewidth{0.95\textwidth}{\rule{\temptablewidth}{1.0pt}}
\begin{tabular*}{\temptablewidth}{@{\extracolsep{\fill}}c|c|c|c|c|c}
method  &  CPU    & GPU (TITAN RTX) & Speedup   & GPU (Tesla V100)  & Speedup\\
\hline
$16^3$   & 28s        &12s   &2.3 &4s&7.0 \\
\hline
$32^3$   & 394s      &82s     &4.8 &39s&10.1 \\
\hline
$64^3$   &5522s        &848s      &6.5 &499s&11.1 \\
\hline
$128^3$  &85556s     &12692s  &6.7 &7927s&10.8\\
\end{tabular*}
{\rule{\temptablewidth}{1.0pt}}
\end{center}
\caption{\label{GPU-CPU-2}Efficiency comparison: the computational
time and speedup for GPU and CPU computation.}
\end{table}

\subsection{Taylor-Green vortex}
The Taylor-Green vortex is a classical problem in fluid dynamics
developed to study vortex dynamics, turbulent transition, turbulent
decay and energy dissipation process \cite{Case-Brachet}. It has
become an excellent case for the evaluation of turbulent flow
simulation methodologies, and been used by many authors for
high-order method validation \cite{Case-Bull,Case-Debonis}. In the
previous study \cite{GKS-high-2}, this problem is provided to test
the performance of WENO-HGKS up to $1024^3$ uniform cells for the
direct numerical simulation (DNS) of nearly incompressible turbulent
flows, and the performance demonstrated the capability of HGKS as a
powerful DNS tool.  In the computation, the DG-HGKS with $\mathbb
P_2$ and $\mathbb P_3$ elements are tested, and the uniform meshes
with $32^3$, $64^3$ and $96^3$ cells are used. As comparison, the
results of  WENO-HGKS are also provided as benchmark, where the
uniform meshes with $64^3$, $128^3$ and $256^3$ cells are used and
fifth-order WENO scheme with linear weights is adopted.

\begin{figure}[!h]
\centering
\includegraphics[width=0.475\textwidth]{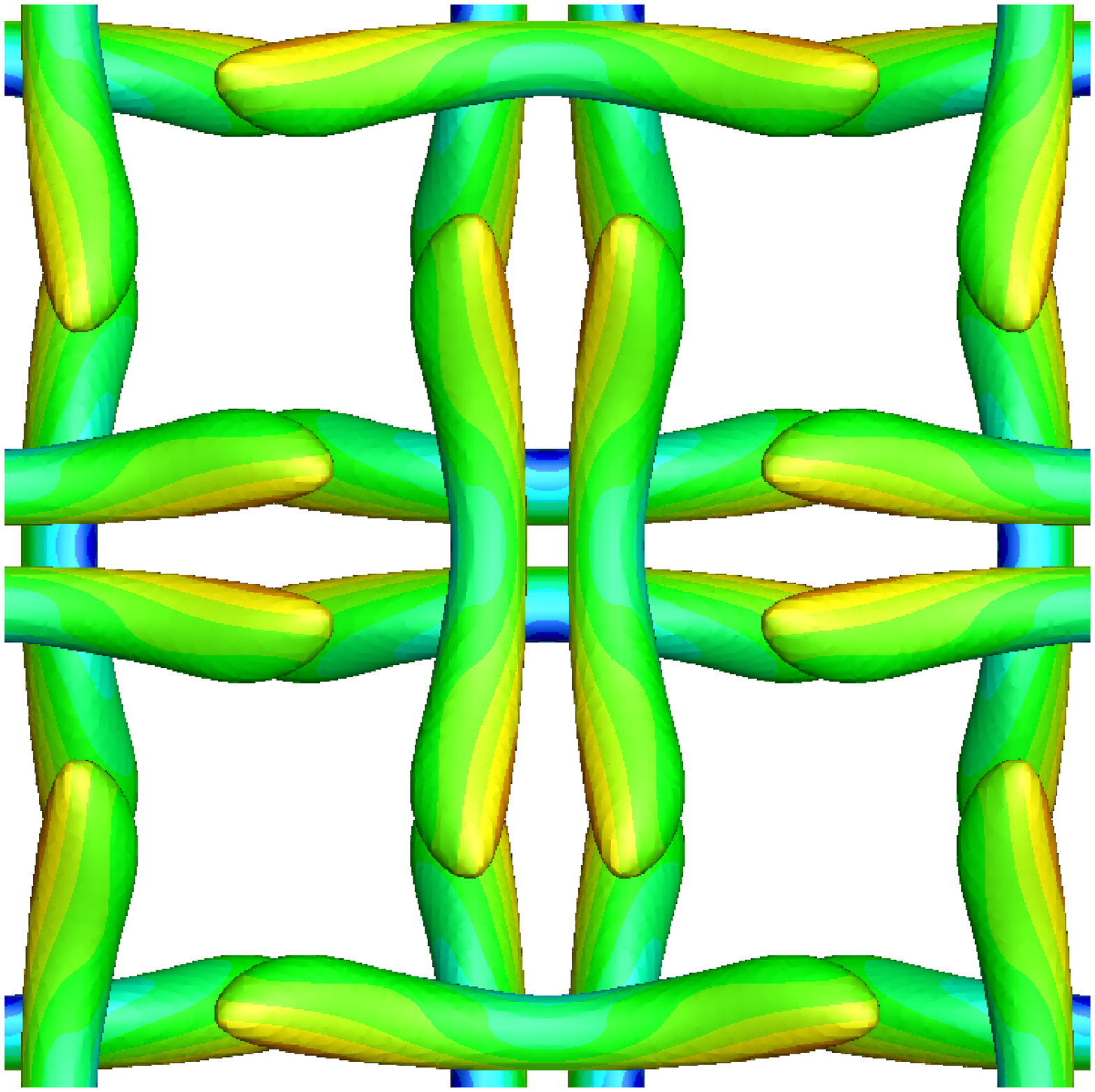}
\includegraphics[width=0.475\textwidth]{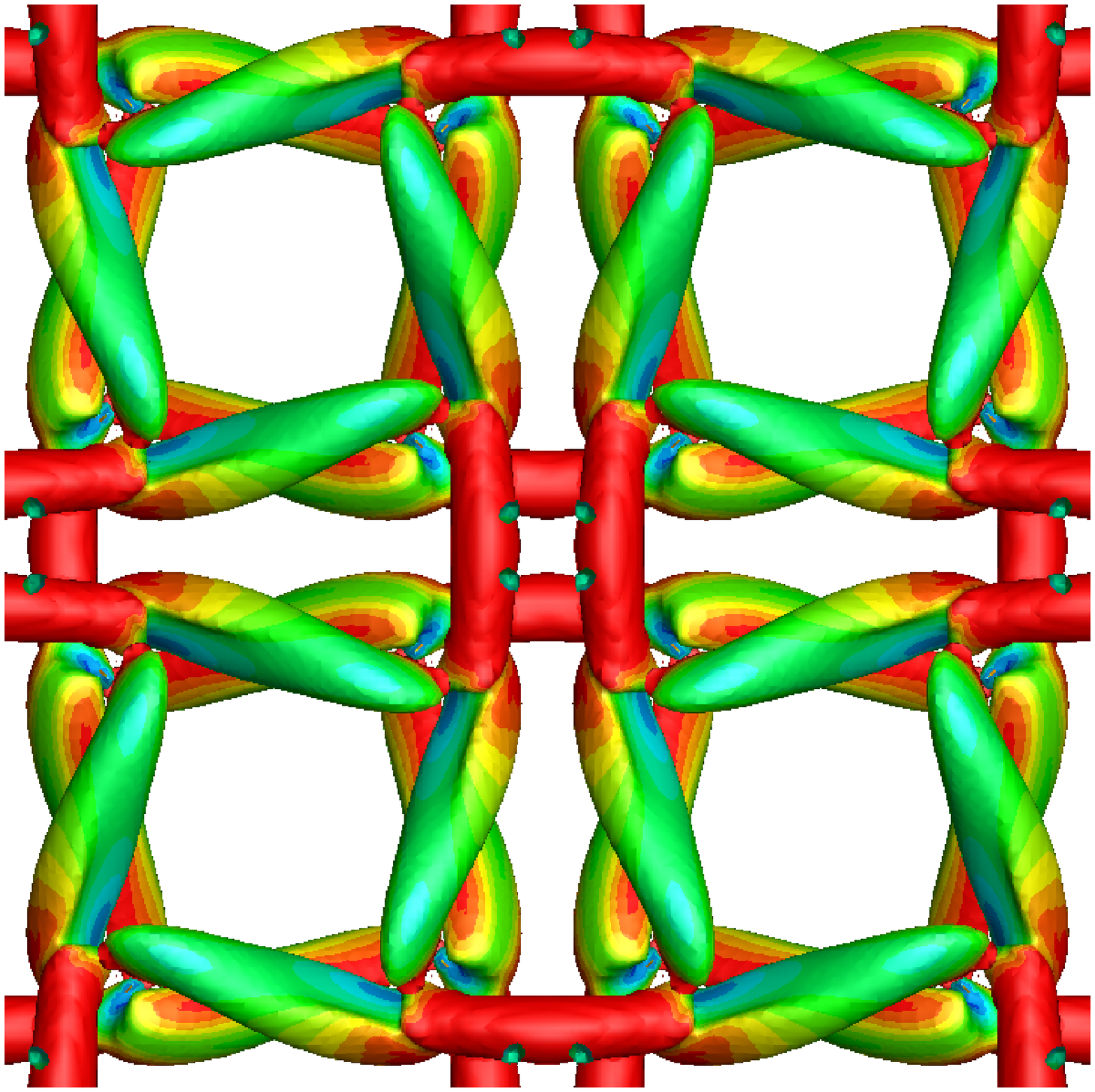}
\includegraphics[width=0.475\textwidth]{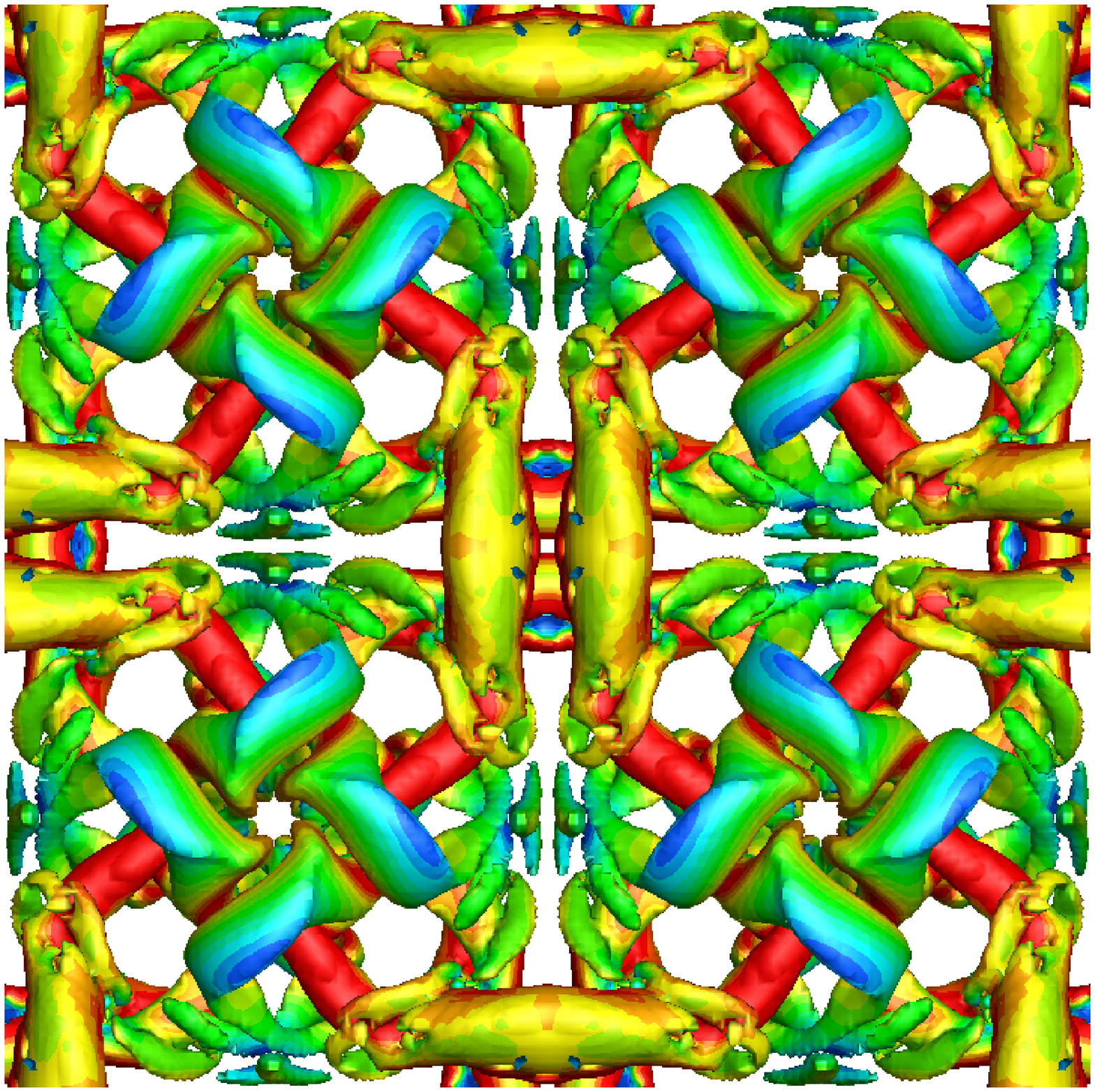}
\includegraphics[width=0.475\textwidth]{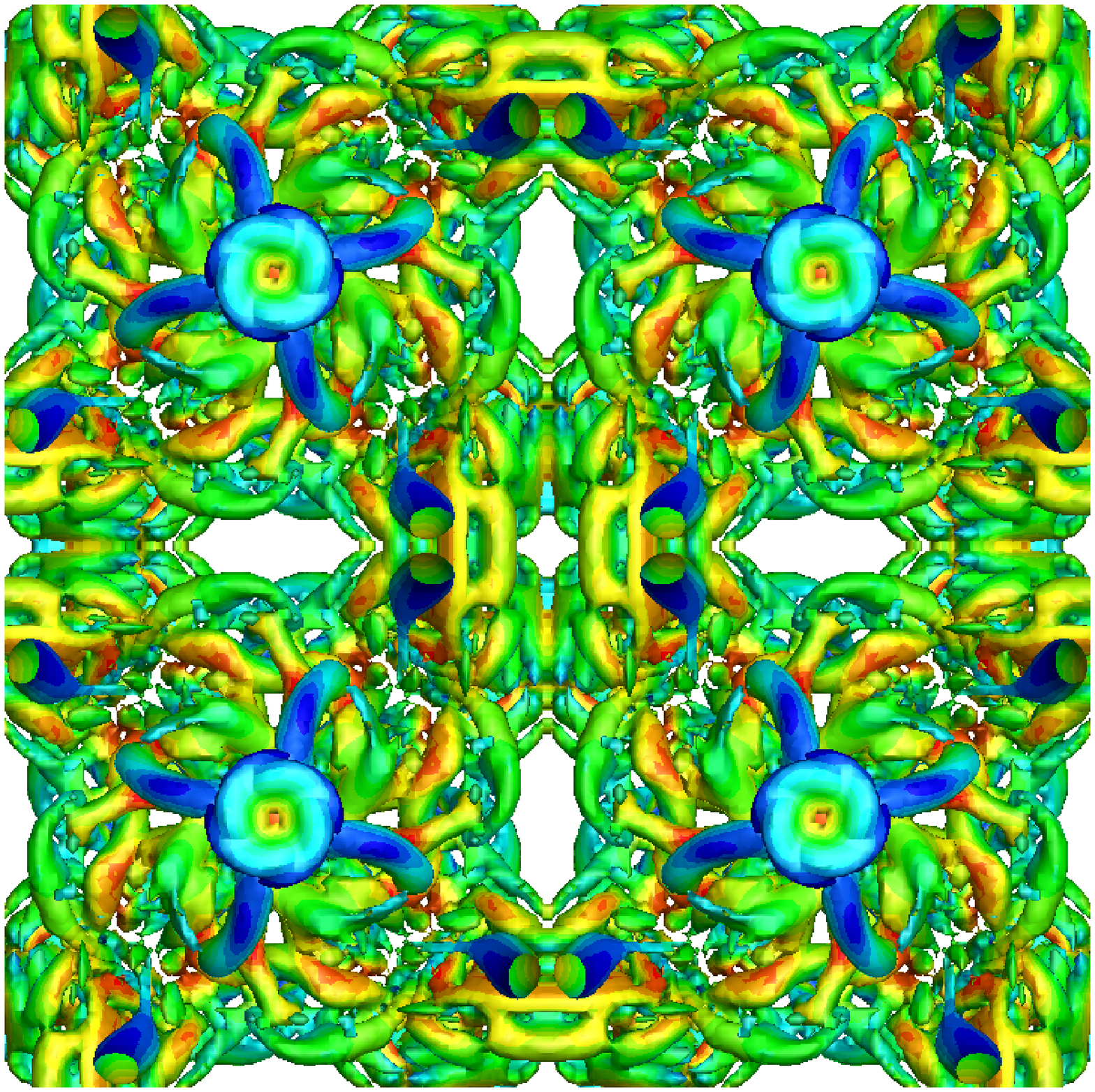}
\caption{\label{Taylor-Green-iso} Taylor-Green vortex: iso-surface
of the second invariant of velocity gradient tensor at $t = 2.5$,
$5$, $7.5$ and $10$ colored by velocity magnitude.}
\end{figure}

The flow is computed within a periodic square box defined as $-\pi
L\leq x, y, z\leq \pi L$. With a uniform temperature, the initial
condition is given by
\begin{align*}
U= & V_0\sin(\frac{x}{L})\cos(\frac{y}{L})\cos(\frac{z}{L}), \\
V= & -V_0\cos(\frac{x}{L})\sin(\frac{y}{L})\cos(\frac{z}{L}), \\
W= & 0, \\
p= & p_0+\frac{\rho_0V_0^2}{16}(\cos(\frac{2x}{L})+\cos(\frac{2y}{L}))(\cos(\frac{2z}{L})+2),
\end{align*}
where $L=1, V_0=1, \rho_0=1$, and the Mach number takes
$M_0=V_0/c_0=0.1$, where $c_0$ is the sound speed. The
characteristic convective time $t_c = L/V_0$, the fluid is a perfect
gas with $\gamma=1.4$, Prandtl number is $Pr=1$ and Reynolds
number $Re=1600$. To test the performance of high-order schemes,
several diagnostic quantities are computed from the flow as it
evolves in time.  This case is given by a simple construction, and
contains several key physical processes including vortex stretching,
interaction and dilatation effects. As time evolves, the vortex
roll-up, stretch and interact, eventually breaking down into
turbulence. The iso-surface of the second invariant of velocity
gradient tensor $Q_v$ colored by velocity magnitude at $t = 2.5, 5,
7.5$ and $10$ with $96^3$ cells for $\mathbb P_3$ element are shown
in Figure.\ref{Taylor-Green-iso}. At the earliest time, the flow
behaves inviscidly as the vortex begin to evolve and roll-up. At
$t=10$, the coherent structures breakdown.

\begin{figure}[!h]
\centering
\includegraphics[width=0.45\textwidth]{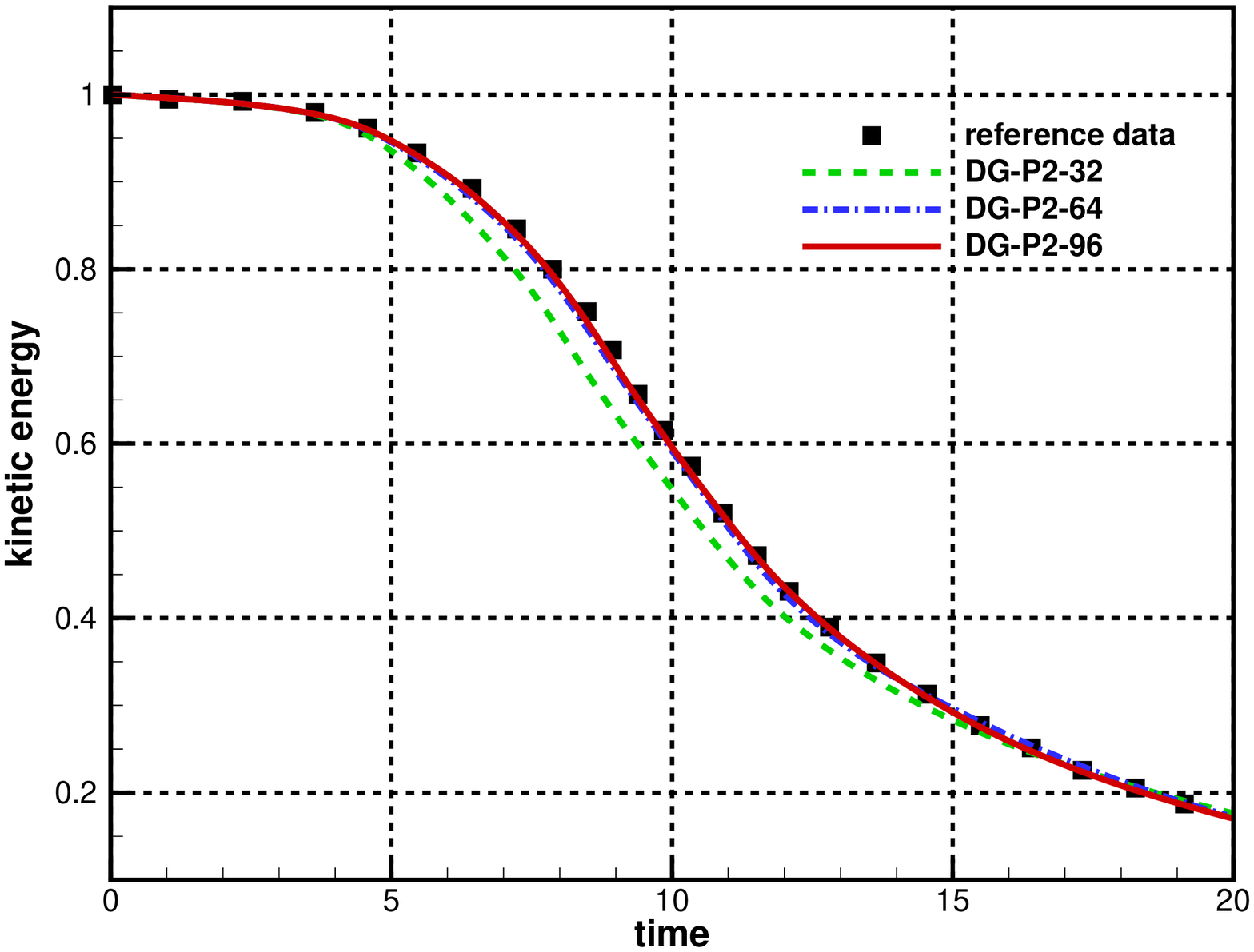}
\includegraphics[width=0.45\textwidth]{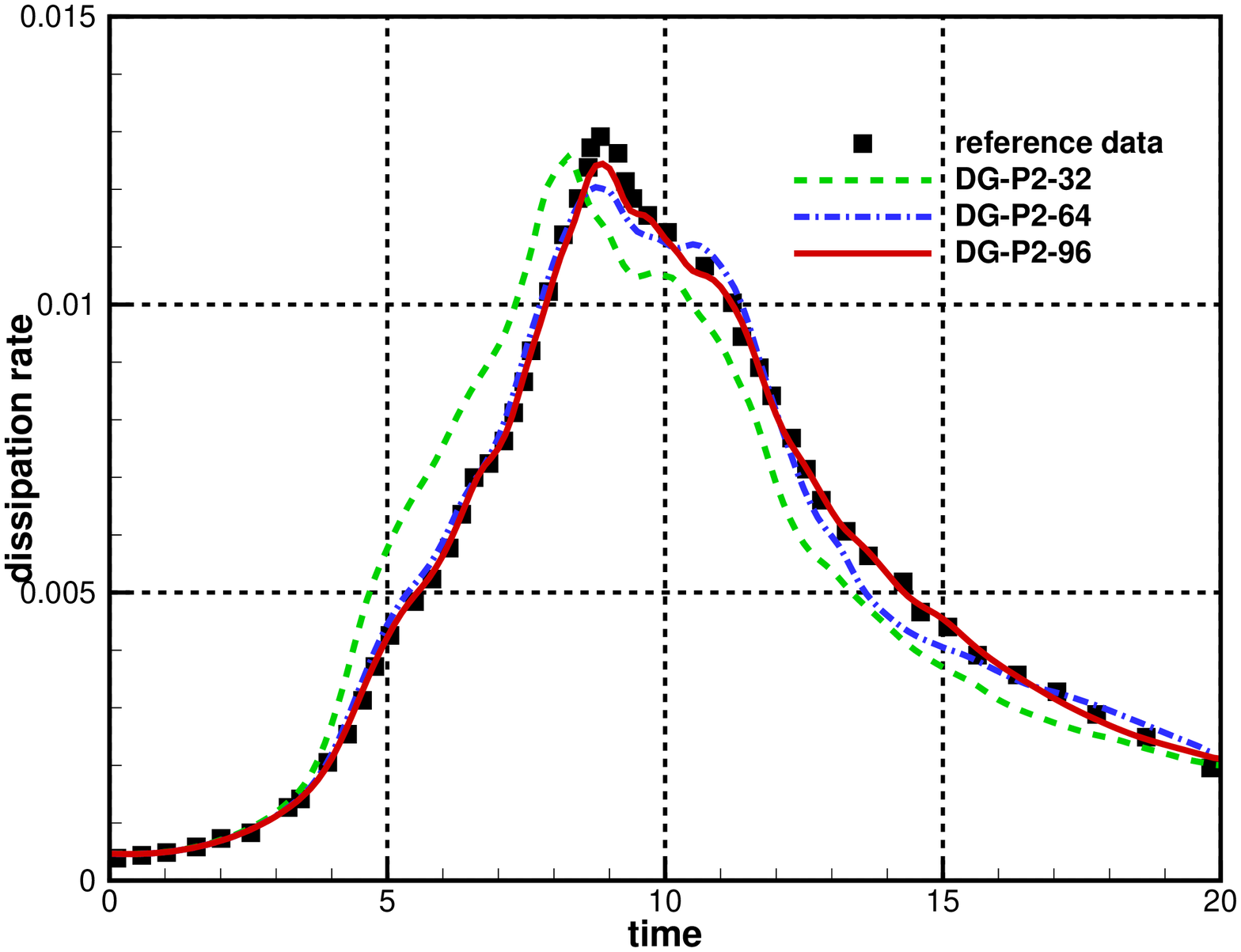}
\caption{\label{Taylor-Green-1A} Taylor-Green vortex: the time
history of kinetic energy (left) and dissipation rate (right) of
DG-HGKS-$\mathbb P_2$  with $32^3$, $64^3$ and $96^3$ cells.}
\includegraphics[width=0.45\textwidth]{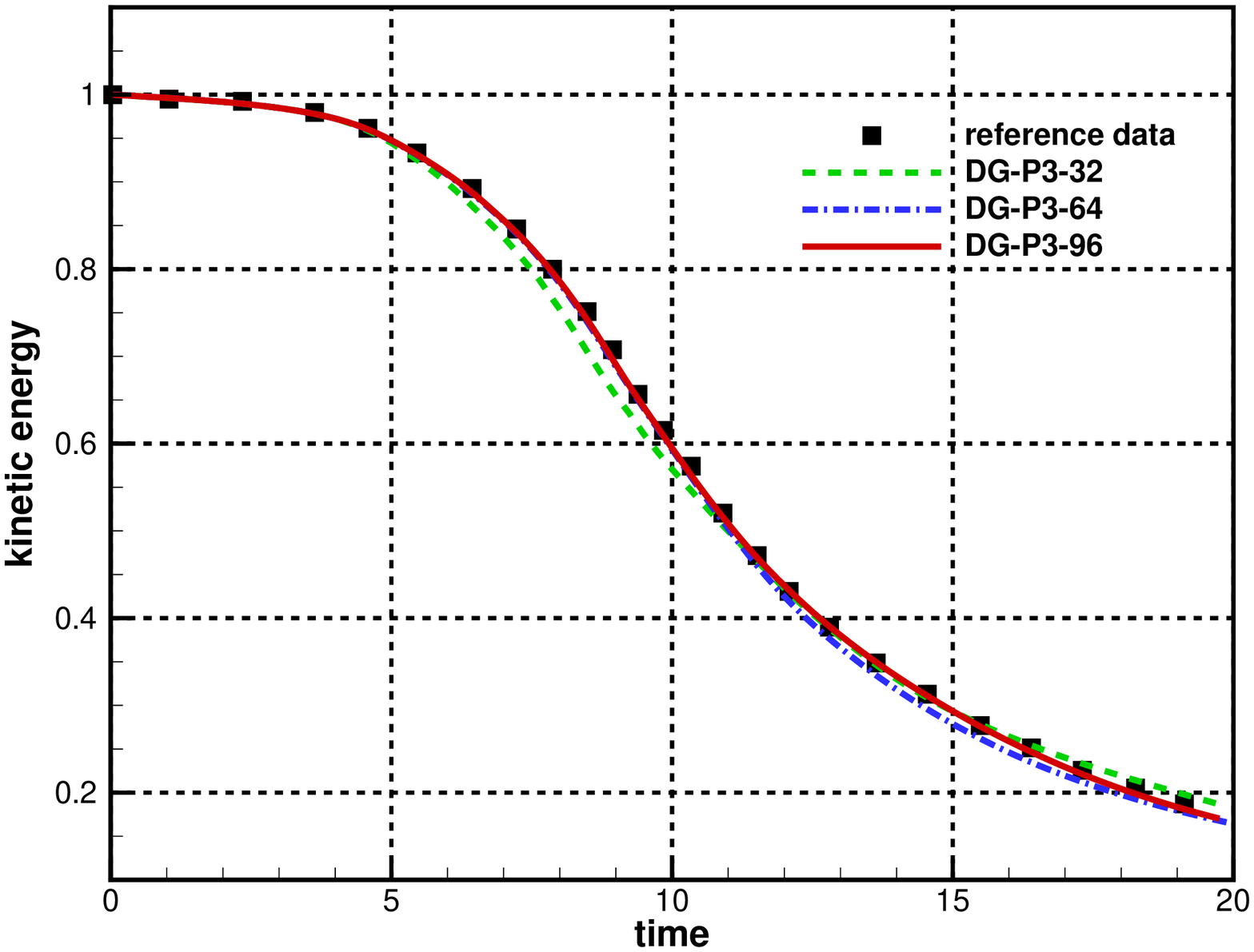}
\includegraphics[width=0.45\textwidth]{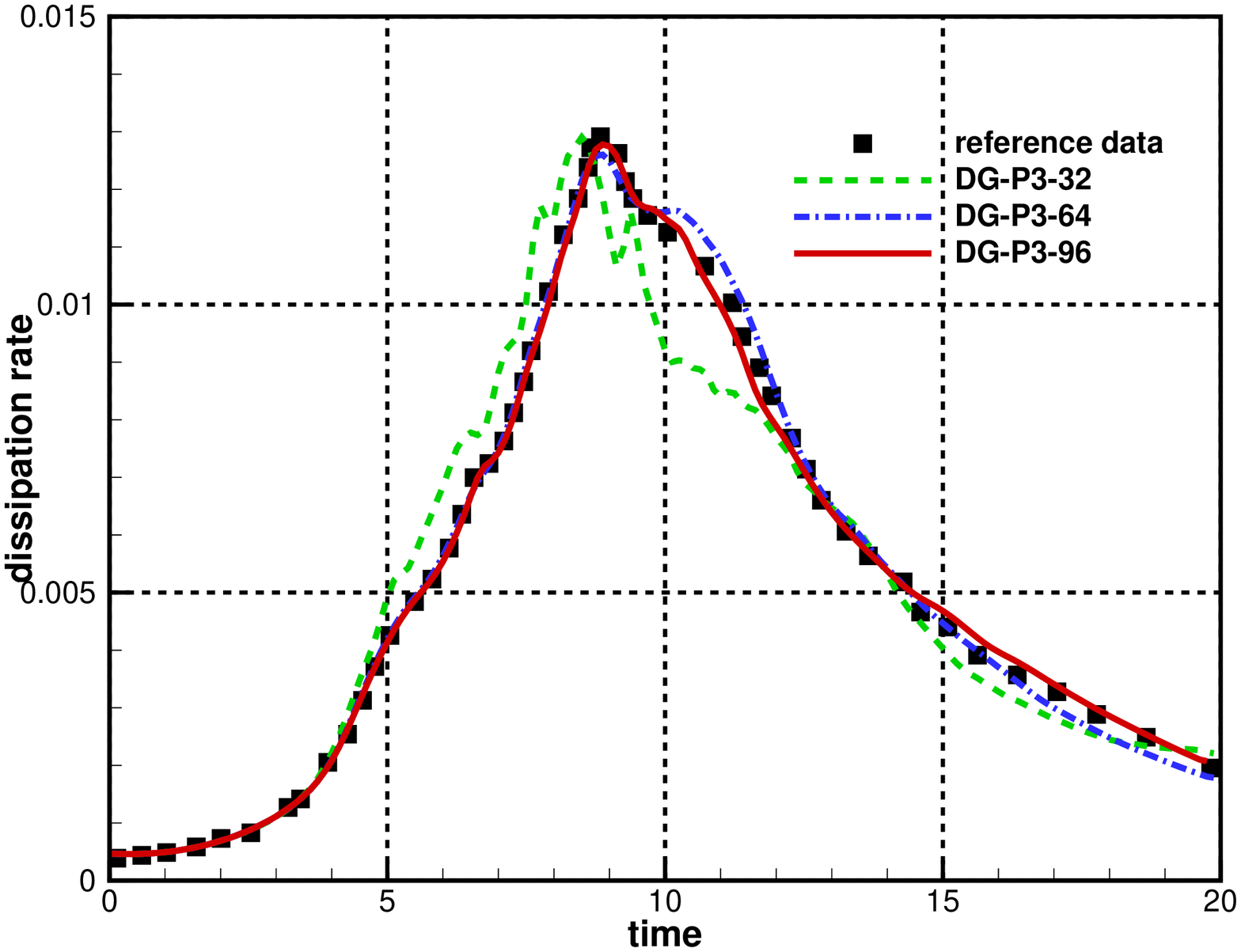}
\caption{\label{Taylor-Green-1B} Taylor-Green vortex: the time
history of kinetic energy (left) and dissipation rate (right) of
DG-HGKS-$\mathbb P_3$ with $32^3$, $64^3$ and $96^3$ cells.}
\includegraphics[width=0.45\textwidth]{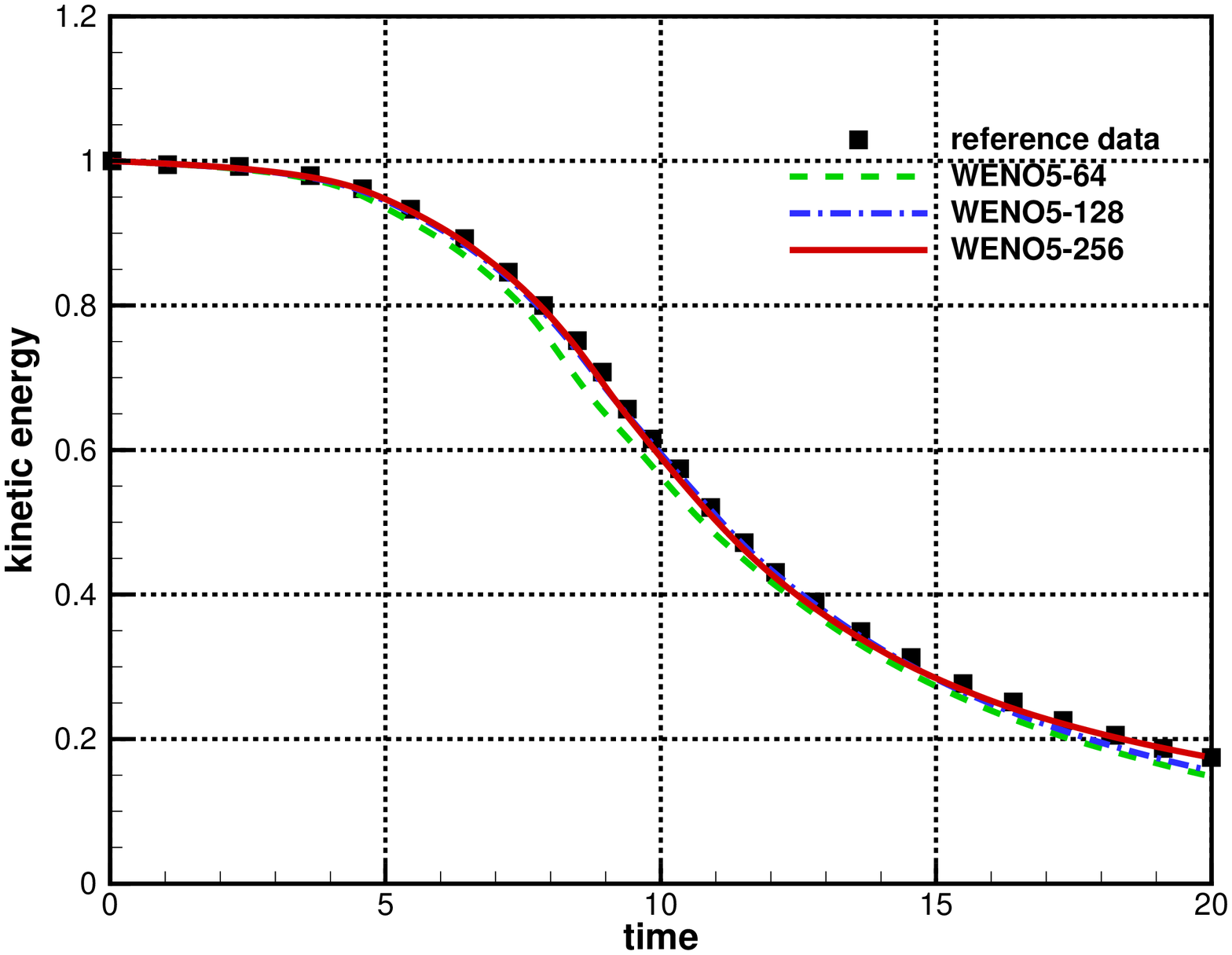}
\includegraphics[width=0.45\textwidth]{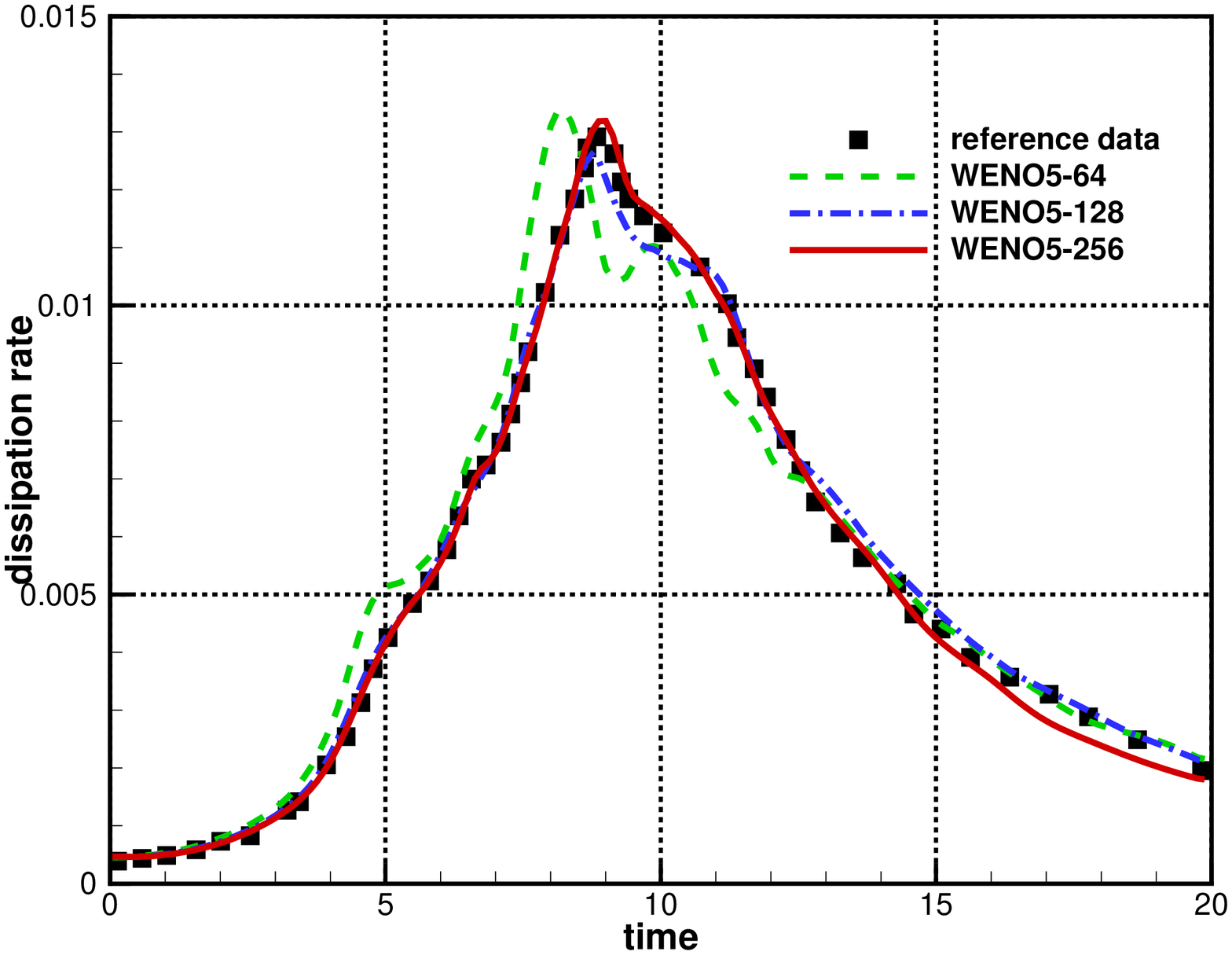}
\caption{\label{Taylor-Green-1C} Taylor-Green vortex: the time
history of kinetic energy (left) and dissipation rate (right) of
WENO-HGKS with $64^3$, $128^3$ and $256^3$ cells.}
\end{figure}

\begin{figure}[!h]
\centering
\includegraphics[width=0.45\textwidth]{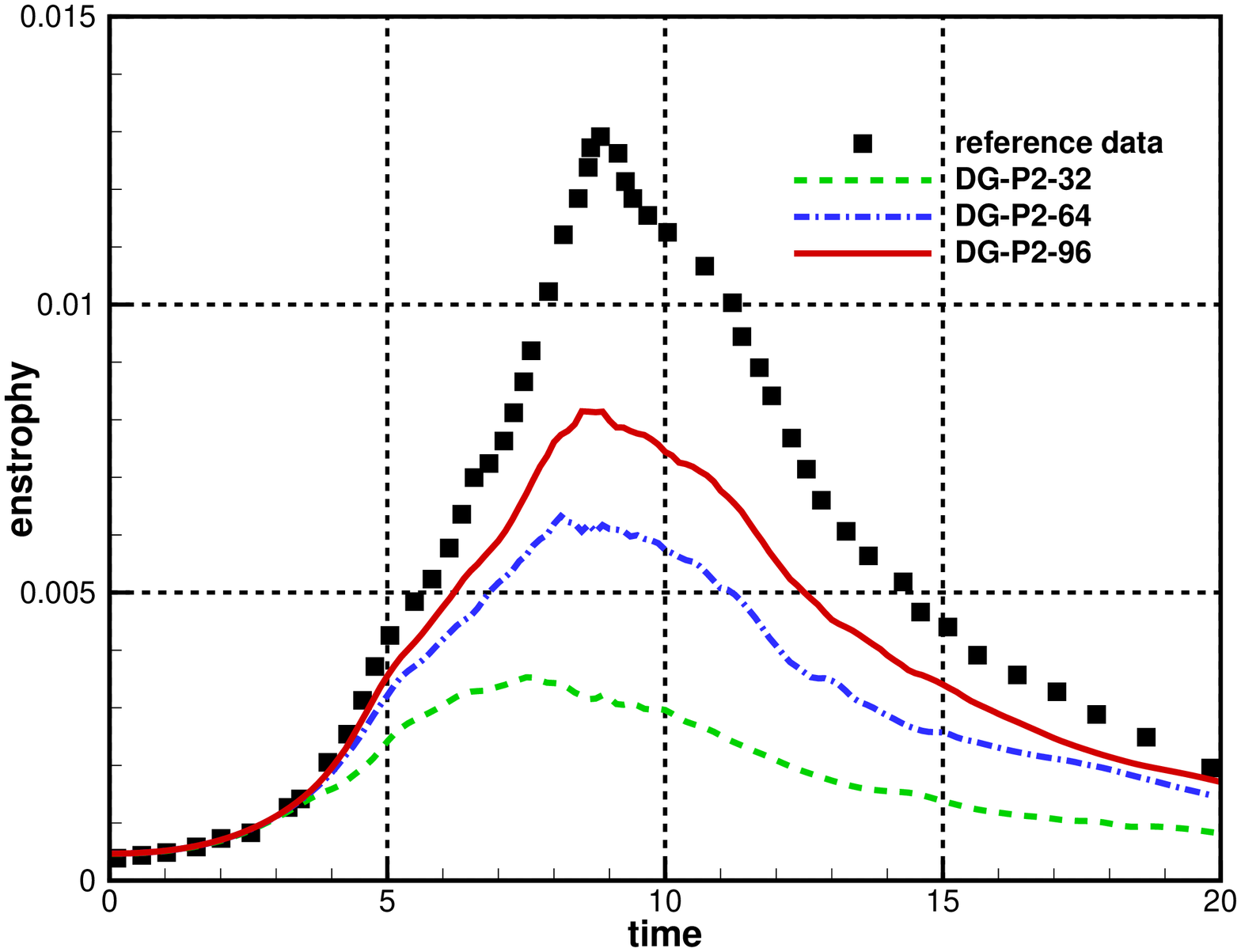}
\includegraphics[width=0.45\textwidth]{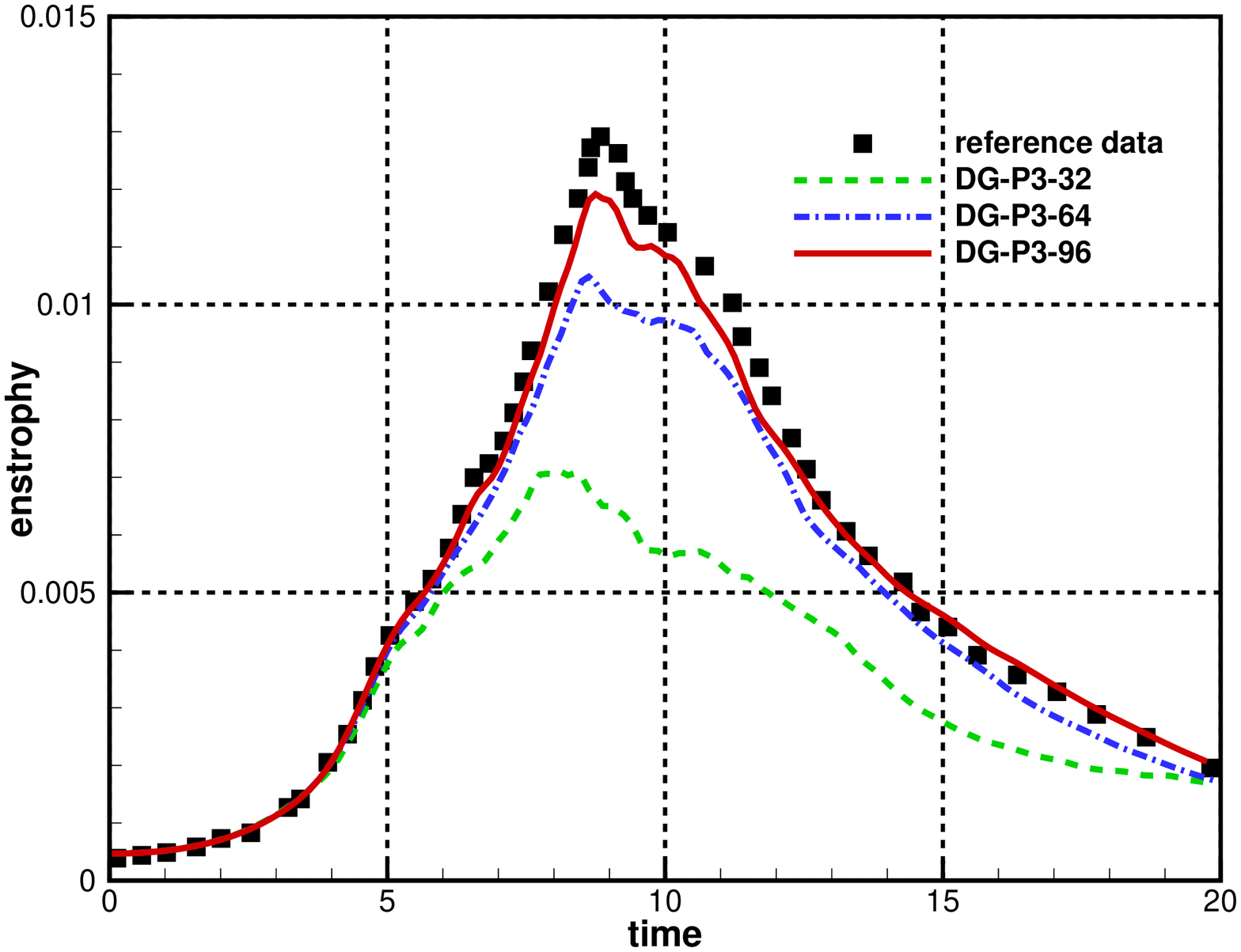}
\caption{\label{Taylor-Green-3A} Taylor-Green vortex: the time
history of integrated enstrophy for DG-HGKS-$\mathbb P_2$ (left) and
DG-HGKS-$\mathbb P_3$ (right) with $32^3$, $64^3$ and $96^3$ cells.}
\includegraphics[width=0.45\textwidth]{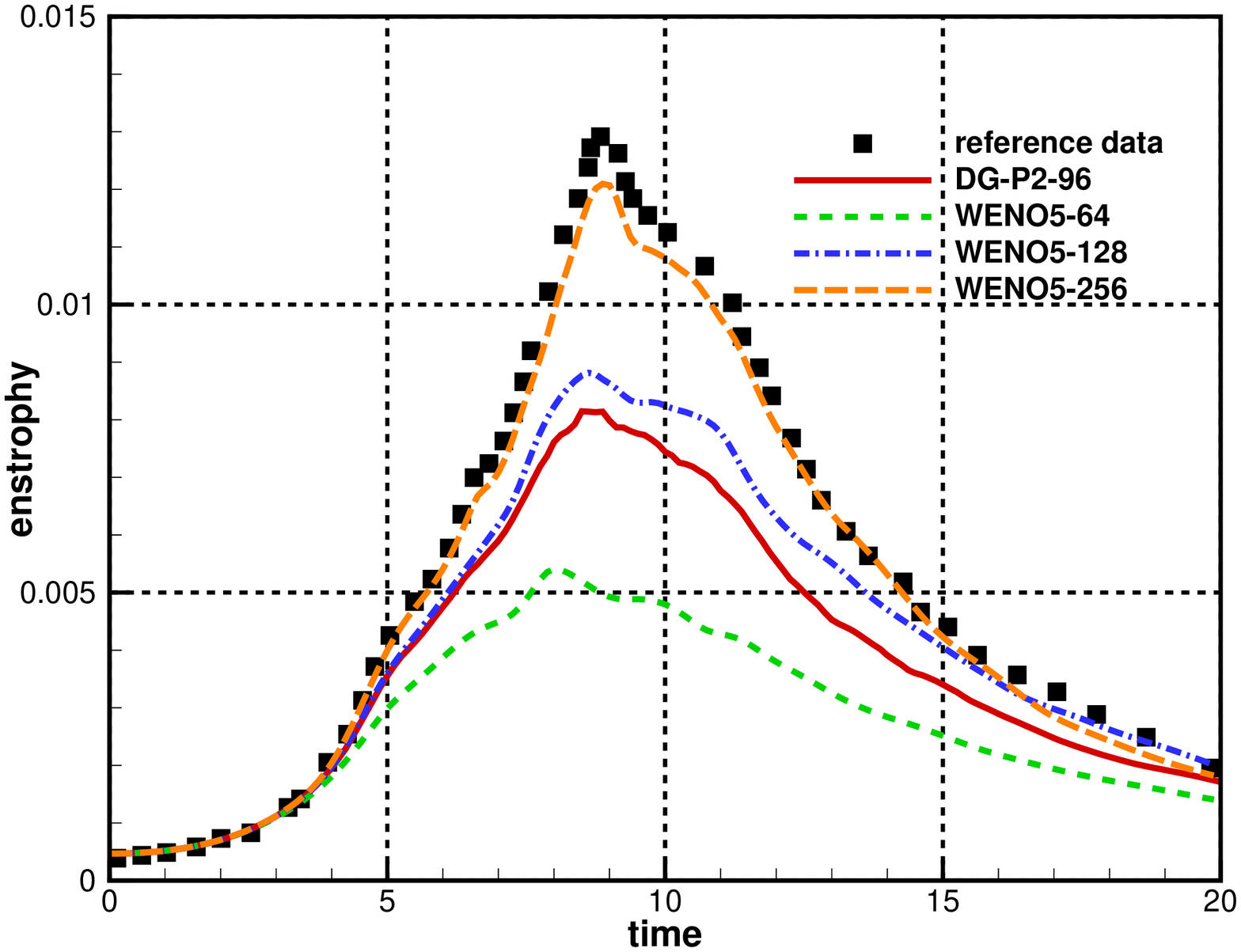}
\includegraphics[width=0.45\textwidth]{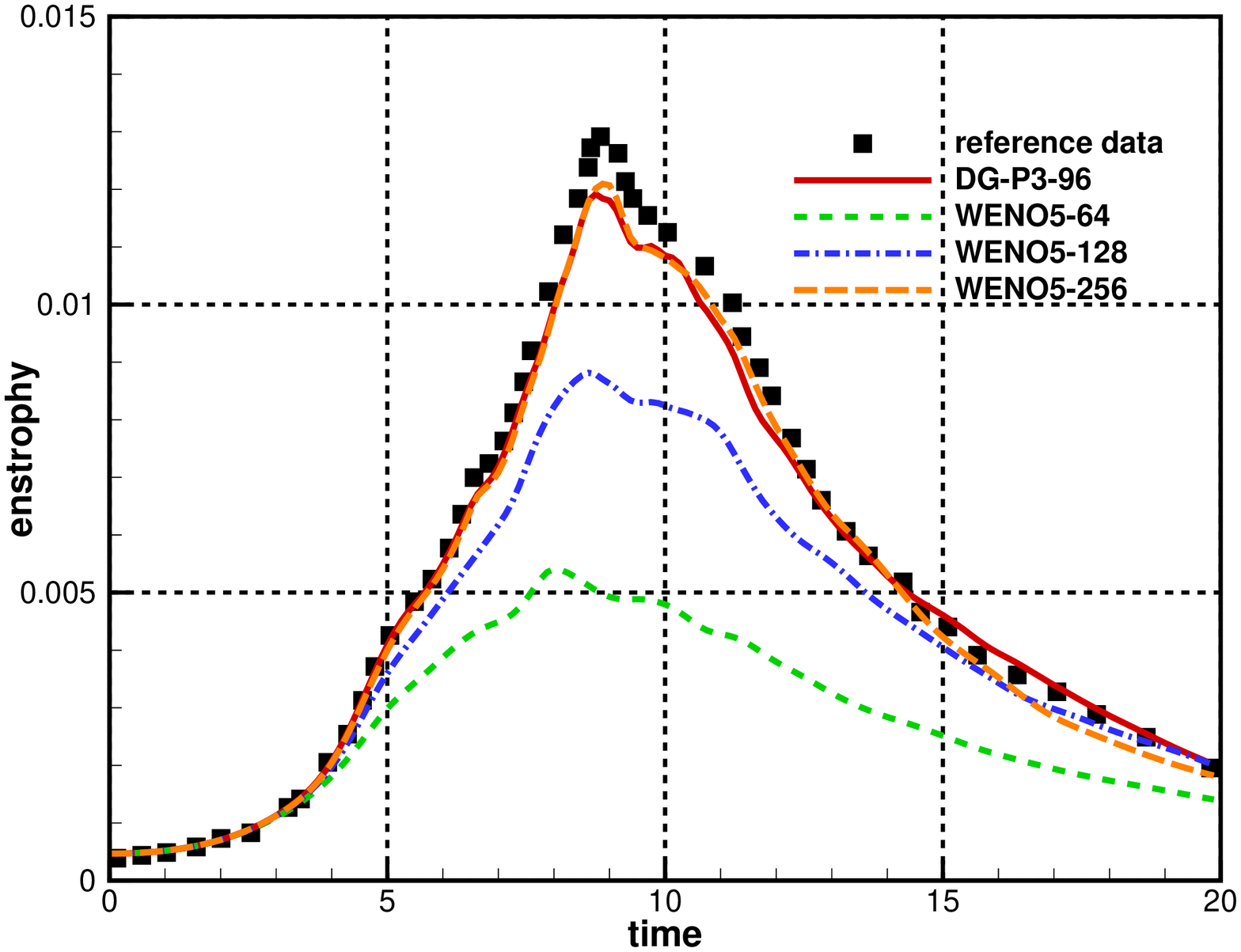}
\caption{\label{Taylor-Green-3B} Taylor-Green vortex: the comparison
of DG-HGKS-$\mathbb P_2$, DG-HGKS-$\mathbb P_3$  and WENO-HGKS for
time history of integrated enstrophy.}
\end{figure}

The volume-averaged kinetic energy is given by
\begin{align*}
E_k=\frac{1}{\rho_0\Omega}\int_\Omega\frac{1}{2}\rho\boldsymbol{U}\cdot\boldsymbol{U} \text{d} \Omega,
\end{align*}
where $\Omega$ is the volume of the computational domain. The
dissipation rate of kinetic energy can be computed by the temporal
derivative of $E_k$
\begin{align*}
\varepsilon(E_k)=-\frac{\text{d}E_k}{\text{d}t},
\end{align*}
which is computed by second order central difference in the
numerical results of $E_k$. For the incompressible limit, the
dissipation rate is related to the integrated enstrophy by
\begin{align*}
\varepsilon(\zeta)=\frac{2\mu}{\rho_0\Omega}\int_\Omega\frac{1}{2}\rho\boldsymbol{\omega}\cdot\boldsymbol{\omega} \text{d} \Omega,
\end{align*}
where $\mu$ is the coefficient of viscosity and $\boldsymbol{\omega}=\nabla\times \boldsymbol{U}$.
The numerical results of BB13 scheme \cite{Case-Debonis,Case-Bogey}
are provided as reference data, where a high-order finite difference
method equipped with fourth-stage third-order algorithm for time
discretization and 13-point stencils for spatial discretization. The
time histories of kinetic energy and dissipation rate of
DG-HGKS-$\mathbb P_2$ and DG-HGKS-$\mathbb P_3$ and WENO methods are
shown in Figure.\ref{Taylor-Green-1A}, Figure.\ref{Taylor-Green-1B}
and Figure.\ref{Taylor-Green-1C}, respectively. The reasonable
agreements are observed with the reference solution, except for the
results of DG schemes with $32^3$ cells and WENO scheme with $64^3$
cells. With the refinement of mesh, the DG-HGKS can capture the
benchmark solution more accurately, which verifies the convergence
of DG schemes. The dissipation rate is determined by both physical
and numerical dissipation, and most of high-order schemes can
resolve it well with a not very refine mesh. However, the integrated
enstrophy is much more difficult to capture for all schemes. The
time histories of  integrated enstrophy are shown in
Figure.\ref{Taylor-Green-3A} for DG-HGKS-$\mathbb P_2$ and
DG-HGKS-$\mathbb P_3$. The magnitude of errors highly depends on the
order of accuracy of the scheme. A large discrepancy is observed for
the integrated enstrophy except for DG-HGKS-$\mathbb P_3$ using
$96^3$ cells. The comparisons with WENO scheme are given in
Figure.\ref{Taylor-Green-3B} for the integrated enstrophy. According
to the quantitative analysis of numerical dissipation
\cite{GKS-high-2}, the performance of the integrated enstrophy is
sensitive to both numerical dissipation and degree of freedom (DoF)
of high-order scheme. The DG-HGKS-$\mathbb P_2$ has $10$ and
DG-HGKS-$\mathbb P_3$ has $20$ DoFs inside each cell. Meanwhile, the
WENO-HGKS has only one DoF inside each cell. Due to the only
third-order of accuracy, the performance of DG-$\mathbb P_2$ method
with $96^3$ cells are only comparable with WENO-HGKS with $128^3$
cells, and a large discrepancy is also observed compared with the
reference data. With the increase of DoF and order of accuracy,
DG-HGKS-$\mathbb P_3$ performs much better than DG-HGKS-$\mathbb
P_2$. Due to the approximate number of DoF,  the time history of
integrated enstrophy of DG-HGKS-$\mathbb P_3$ with $96^3$ cells are
comparable with that of WENO-HGKS with $256^3$ cells, and both of
them agree well with the reference data.
However, for the flow with discontinuities, such as the compressible isotropic turbulence with
initial turbulent Mach number $0.5$ and Taylor microscale Reynold
number $72$, WENO-HGKS works well even with linear weights, but the
DG-HGKS fails without special treatment. The DG methods seem to lack
robustness, and more great efforts need to be paid to the limiters
and trouble cell indicators beforehand
\cite{DG-limiter-1,DG-limiter-2}.

The comparison of computational cost is also studied, where the
detailed mesh size, time step, CFL number and CPU time for DG-HGKS
and WENO-HGKS are presented in Table.\ref{time_table2}. The CPU time
are computed per 10 steps for each method by Intel i7-9700 CPU using
OpenMP directives. Unsurprisingly, due to the fewer Gaussian quadrature points of 
DG-HGKS-$\mathbb P_2$, DG-HGKS-$\mathbb P_2$ is much more efficient 
DG-HGKS-$\mathbb P_3$. With the same number of cells, no doubt that DG-HGKS is less efficient
than WENO-HGKS. Taken the degree of freedom into account, the CPU
time of WENO-HGKS with $256^3$ cells is around $2.15$ times of
DG-HGKS-$\mathbb P_3$ with $96^3$ cells. The CPU time of WENO-HGKS
is around $1/1.66$ of DG-HGKS-$\mathbb P_3$. Thus, the total
computational cost of WENO-HGKS is comparable with DG-HGKS-$\mathbb
P_3$. However, to improve the accuracy of DG scheme, much more
quadrature points and numerical fluxes are needed,  and much smaller
time step is used, which cause the dramatical increase of
computational cost for DG method.

\begin{table}[!h]
\begin{center}
\def\temptablewidth{0.85\textwidth}{\rule{\temptablewidth}{1.0pt}}
\begin{tabular*}{\temptablewidth}{@{\extracolsep{\fill}}c|c|c|c|c}
method  &   mesh size & time    step   &  CFL number     & CPU time (per 10 steps)\\
\hline
DG-$\mathbb P_2$    &  $96^3$    & $8.86\times10^{-4}$   &  0.15    &  46s\\
\hline
DG-$\mathbb P_3$    &  $96^3$    & $5.35\times10^{-4}$  &  0.09   &  155s\\
\hline
WENO-$5$                &   $256^3$  & $8.92\times10^{-4}$   & 0.40  &    333s  \\
\end{tabular*}
{\rule{\temptablewidth}{1.0pt}}
\end{center}
\caption{\label{time_table2}  Taylor-Green vortex: comparison of
efficiency for DG-HGKS and WENO-HGKS.}
\end{table}

\section{Conclusion}
In this paper, the DG-HGKS is developed for the compressible Euler
and Navier-Stokes equations. Different from the traditional DG
approaches, the time dependent flux is provided by the gas-kinetic
flow solver. In the framework of DG-HGKS, the inviscid and viscous
fluxes can be calculated uniformly, which avoid the treatment of
diffusion terms in LDG \cite{LDG} and DDG \cite{DDG2} methods. Due
to the time dependent numerical flux, the temporal accuracy is
achieved by two-stage fourth-order discretization
\cite{GRP-high-1,GKS-high-1}. With $\mathbb P_2$ and $\mathbb P_3$
elements, the third-order and fourth-order spatial accuracy can be
achieved. In the computation, we mainly concentrate on the smooth
flows, and numerical tests are presented to validate the performance
of current scheme.  The accuracy tests are presented to validate the
optimal convergence and super-convergence for DG discrezation with
uniform and nonuniform meshes, and the Taylor-Green vortex are
provided to validate the efficiency and accuracy for the nearly
incompressible turbulence. The numerical results of WENO-HGKS are
given as well. Both DG-HGKS and WENO-HGKS work well for the smooth
flows. As expected, the numerical performances, including the
resolution and efficiency, are comparable with the approximate
number of degree of freedom. To accelerate the computation, the
DG-HGKS is implemented to run on Nvidia GPUs using CUDA. Obtained
results are compared with those obtained by Intel i7-9700 CPU using
OpenMP directives. The GPU code achieves 6x-7x speedup with TITAN
RTX, and 10x-11x speedup with Tesla V100. In the future, more
challenging problems for compressible flows, such as the supersonic
and hypersonic flat plate turbulent boundary layer, will be
investigated.

\section*{Ackonwledgement}
This research is supported by National Natural Science Foundation of
China (11701038), the Fundamental Research Funds for the Central
Universities.


\begin{thebibliography}{10}

\bibitem{BR1} F. Bassi, A. Crivellini, D.A. Di Pietro, S. Rebay, An artificial compressibility flux for the discontinuous Galerkin solution of the incompressible Navier-Stokes equations, J. Comput. Phys. 218 (2006) 794-815.

\bibitem{BR2} F. Bassi, A. Crivellini, D.A. Di Pietro, S. Rebay, An implicit high-order discontinuous Galerkin method for steady and unsteady incompressible flows, Comput. Fluids 36 (2007) 1529-1546.

\bibitem{BGK-1} P.L. Bhatnagar, E.P. Gross, M. Krook, A Model for Collision Processes in Gases I: Small Amplitude Processes in Charged and Neutral One-Component Systems, Phys. Rev. 94 (1954) 511-525.

\bibitem{Case-Bogey} C. Bogey, C. Bailly, A family of low dispersive and low dissipative explicit schemes for flow and noise computations, J. Comput. Phys. 194 (2004) 194-214.

\bibitem{Case-Brachet} M.E. Brachet, D.I. Meiron, S.A. Orszag, B.G. Nickel, R.H. Morf, U. Frisch, Small-scale structure of the Taylor-Green vortex, J. Fluid. Mech. 130 (1983) 411-452.

\bibitem{Case-Bull} J. R. Bull, A. Jameson, Simulation of the compressible Taylor-Green vortex using high-order flux reconstruction schemes, AIAA 2014-3210.

\bibitem{WENO-Z} R. Borges, M. Carmona, B. Costa, W. S. Don, An improved weighted essentially non-oscillatory scheme for hyperbolic conservation laws, J. Comput. Phys. 227 (2008) 3191-3211.

\bibitem{GKS-high-4}  G.Y. Cao, L. Pan, K. Xu, High-order gas-kinetic scheme with parallel computation for direct numerical simulation of turbulent flows, J. Comput. Phys.  448 (2022) 110739

\bibitem{DG-superconvergence}  W. Cao, C.-W. Shu, Y. Yang and Z. Zhang, Superconvergence of discontinuous Galerkin method for scalar nonlinear hyperbolic equations. SIAM J. Numer. Anal. 56 (2018) 732-765.

\bibitem{BGK-2} S. Chapman, T.G. Cowling, The Mathematical theory of non-uniform gases, third edition, Cambridge University Press, (1990).

\bibitem{DDG2}  J. Cheng, X.Q. Yang, X.D. Liu, T.G. Liu, H. Luo, A direct discontinuous Galerkin method for the compressible Navier–Stokes equations on arbitrary grids, J. Comput. Phys. 327 (2016) 484-502.

\bibitem{DG1}  B. Cockburn, C.-W. Shu, TVB Runge-Kutta local projection discontinuous Galerkin finite element method for conservation laws II: General framework, Math. Comput. 52 (1989) 411-435.

\bibitem{DG3}  B. Cockburn, C.-W. Shu, The Runge-Kutta discontinuous Galerkin method for conservation laws V: Multidimensional systems, J. Comput. Phys. 141 (1998) 199-224.

\bibitem{LDG}  B. Cockburn, G. Kanschat, D. Sch\"{o}tzau, A locally conservative LDG method for the incompressible Navier-Stokes equations, Math. Comput. 74 (2005) 1067-1095.

\bibitem{DG2}  B. Cockburn, S.Y. Lin, C.-W. Shu, TVB Runge-Kutta local projection discontinuous Galerkin finite element method for conservation laws III: One dimensional systems, J. Comput. Phys. 84 (1989) 90-113.

\bibitem{Case-Debonis} J. Debonis, Solutions of the Taylor-Green vortex problem using high-resolution explicit finite difference methods, AIAA Paper (2013) 2013-0382.

\bibitem{GRP-high-2} Z.F. Du, J.Q. Li, A Hermite WENO reconstruction for fourth order temporal accurate schemes based on the GRP solver for hyperbolic conservation laws, J. Comput. Phys. 355 (2018) 385-396.

\bibitem{CPR-1} H. Gao, Z.J. Wang, A high-order lifting collocation penalty formulation for the Navier-Stokes equations on 2D mixed grids, AIAA Paper (2009) 3784.

\bibitem{TVD-RK} S. Gottlieb, C. W. Shu, Total variation diminishing Runge-Kutta schemes, Mathematics of computation, 67 (1998) 73-85.

\bibitem{ENO-1} A. Harten, B. Engquist, S. Osher and S. R. Chakravarthy, Uniformly high order accurate essentially non-oscillatory schemes, III. J. Comput. Phys. 71 (1987) 231-303.

\bibitem{FR} H.T. Huynh, A flux reconstruction approach to high-order schemes including discontinuous Galerkin methods, AIAA Paper 2007 4079.
   
\bibitem{CPR-2}  H.T. Huynh, Z.J. Wang, P.E. Vincent,  High-order methods for computational fluid dynamics: A brief review of compact differential formulations on unstructured grids, Computers $\&$ Fluids 98 (2014) 209-220

\bibitem{GKS-high-3} X. Ji,  L. Pan, W. Shyy , K. Xu,  A compact fourth-order gas-kinetic scheme for the Euler and Navier-Stokes equations, J. Comput. Phys. 372 (2018) 446-472.

\bibitem{GKS-high-2} X. Ji, F.X. Zhao, W. Shyy, K. Xu,   A family of high-order gas-kinetic schemes and its comparison with Riemann solver based high-order methods,  J. Comput. Phys. 356 (2018) 150-173.

\bibitem{WENO-JS} G.S. Jiang, C.-W. Shu, Efficient implementation of weighted ENO schemes, J. Comput. Phys. 126 (1996) 202-228.

\bibitem{GRP-high-1} J.Q. Li, Z.F. Du, A two-stage fourth order time-accurate discretization for Lax-Wendroff type flow solvers I. hyperbolic conservation laws, SIAM J. Sci. Computing, 38 (2016) 3046-3069.

\bibitem{DDG1}  H.L. Liu, J. Yan, The direct discontinuous Galerkin (DDG) methods for diffusion problems, SIAM J. Numer. Anal. 47 (2009) 675-698.

\bibitem{WENO-Liu} X.D. Liu, S. Osher, T. Chan, Weighted essentially non-oscillatory schemes, J. Comput. Phys. 115 (1994) 200-212.

\bibitem{SD} Y. Liu, M. Vinokur, Z. J. Wang, Spectral difference method for unstructured grids I: Basic formulation, J. Comput. Phys. 216 (2006) 780-801.

\bibitem{RDG}  H. Luo, L.Q. Luo, R. Nourgaliev, V.A. Mousseau, N. Dinh, A reconstructed discontinuous Galerkin method for the compressible Navier-Stokes equations on arbitrary grids, J. Comput. Phys. 229 (2010) 6961-6978.

\bibitem{DG-convergence} X. Meng, C.-W. Shu, B. Wu, Optimal error estimates for discontinuous Galerkin methods based on upwind-biased fluxes for linear hyperbolic equations. Math. Comput. 85 (2016) 1225-1261.

\bibitem{GKS-high-1} L. Pan, K. Xu, Q.B. Li, J.Q. Li, An efficient and accurate two-stage fourth-order gas-kinetic scheme for the Navier-Stokes equations, J. Comput. Phys. 326 (2016) 197-221.

\bibitem{DG-limiter-1}  J. X. Qiu, C.-W. Shu, Hermite WENO schemes and their application as limiters for Runge–Kutta discontinuous Galerkin method: one dimensional case, J. Comput. Phys.  193 (2003) 115-135.

\bibitem{DG-limiter-2}  J.X. Qiu, C.-W. Shu, Runge–Kutta discontinuous Galerkin method using WENO limiters, SIAM Journal on Scientific Computing 26 (2005) 907-929.

\bibitem{DG0}  W.H. Reed, T.R. Hill, Triangular mesh methods for the neutron transport equation, Los Alamos Scientific Laboratory Report, LA-UR-73-479 (1973).

\bibitem{GKS-DG-2} X.D. Ren, K. Xu, W. Shyy, C.W. Gu,  A multi-dimensional high-order discontinuous Galerkin method based on gas kinetic theory for viscous flow computations, J. Comput. Phys.  292 (2015) 176-193.

\bibitem{ENO-2} C.-W. Shu, S. Osher, Efficient implementation of essentially non-oscillatory shock capturing schemes, J. Comput. Phys. 77 (1988) 439-471.

\bibitem{Riemann-appro} E.F. Toro, Riemann Solvers and Numerical Methods for Fluid Dynamics, Third Edition, Springer (2009).

\bibitem{SV} Z.J. Wang, Spectral (finite) volume method for conservation laws on unstructured grids: basic formulation, J. Comput. Phys. 178 (2002) 210-251.

\bibitem{GPU} C.F. Xu, X.G. Deng, et el., Collaborating CPU and GPU for large-scale high-order CFD simulations with complex grids on the TianHe-1A supercomputer, J. Comput. Phys. 278 (2014) 275-297.

\bibitem{GKS-Xu1} K. Xu, Direct modeling for computational fluid dynamics: construction and application of unfied gas kinetic schemes, World Scientific (2015).

\bibitem{GKS-Xu2} K. Xu, A gas-kinetic BGK scheme for the Navier-Stokes equations and its connection with artificial dissipation and Godunov method, J. Comput. Phys. 171 (2001) 289-335.

\bibitem{GKS-DG-1} K. Xu, Discontinuous Galerkin BGK method for viscous flow equations: One-dimensional systems. Siam J. Sci. Comput. 25 (2004) 1941-1963.


\end{thebibliography}
\end{document}